\documentclass[leqno,12pt]{article}
\usepackage{amsmath, amssymb}
\usepackage{amsthm} 
\usepackage{latexsym}
\usepackage{url}
\theoremstyle{plain} 
\newtheorem{theorem}{\indent\sc Theorem}[section]
\newtheorem{lemma}[theorem]{\indent\sc Lemma}
\newtheorem{corollary}[theorem]{\indent\sc Corollary}
\newtheorem{proposition}[theorem]{\indent\sc Proposition}

\theoremstyle{definition} 
\newtheorem{definition}[theorem]{\indent\sc Definition}
\newtheorem{remark}[theorem]{\indent\sc Remark}
\newtheorem{example}[theorem]{\indent\sc Example}

\def\prfb{\begin{proof}}\def\prfe{\end{proof}}
\def\prfofb#1{\begin{proof}[{\sc P{\scriptsize ROOF OF} #1}.]}\def\prfofe{\end{proof}}

\def\prpb{\begin{proposition}}\def\prpe{\end{proposition}}
\def\lemb{\begin{lemma}}\def\leme{\end{lemma}}
\def\thmb{\begin{theorem}}\def\thme{\end{theorem}}
\def\corb{\begin{corollary}}\def\core{\end{corollary}}
\def\xmpb{\begin{example}}\def\xmpe{\end{example}}
\def\dfnb{\begin{definition}}\def\dfne{\end{definition}}
\def\remb{\begin{remark}}\def\reme{\end{remark}}

\def\prpa#1{\label{p:#1}}\def\prpu#1{Proposition~\ref{p:#1}}
\def\lema#1{\label{l:#1}}\def\lemu#1{Lemma~\ref{l:#1}}
\def\thma#1{\label{t:#1}}\def\thmu#1{Theorem~\ref{t:#1}}
\def\cora#1{\label{c:#1}}\def\coru#1{Corollary~\ref{c:#1}}
\def\dfna#1{\label{df:#1}}\def\dfnu#1{\mbox{Definition~\ref{df:#1}}}
\def\seca#1{\label{s:#1}}\def\secu#1{\S~\ref{s:#1}}
\def\eqna#1{\label{e:#1}}\def\eqnu#1{(\ref{e:#1})}

\def\labelenumi{(\roman{enumi})} 
\def\itmb{\begin{enumerate}}\def\itme{\end{enumerate}}
\def\itdb{\begin{itemize}}\def\itde{\end{itemize}}
\def\ittb{\begin{description}}\def\itte{\end{description}}
\def\arrb#1{\begin{array}{#1}}\def\arre{\end{array}}
\def\tabb#1{\par\noindent\begin{tabular}{#1}}
\def\tabe{\end{tabular}\par\noindent}

\def\QED{}
\def\DDD{}

\def\Rom#1{\uppercase\expandafter{\romannumeral#1}}
\def\dsp{\displaystyle}
\def\eps{\epsilon}

\def\limf#1{\displaystyle \lim_{#1\to\infty}}

\def\Ccomb#1#2{\setbox0=\hbox{$\displaystyle\mathrm{C}$}\setbox1=\hbox{%
$\scriptstyle #1$}\kern \wd1{\mathrm{C}}_{\kern -1.05\wd0\kern -0.99\wd1{#1}
 \kern 1.15\wd0{#2}}}

\def\clvec#1#2#3{\def\clvecone{#3}\left(\arrb{c} \dsp #1\\ \dsp #2
 \ifx\clvecone\empty\else\\ \dsp #3\fi\arre\right)}

\def\diff#1#2{\dsp\frac{d\,#1}{d#2}}

\def\le{\leqq} \def\ge{\geqq} 

\def\reals{{\mathbb R}}
\def\nreals#1{{\mathbb R}^{#1}}
\def\preals{[0,\infty)}  \def\preals{{\mathbb R_+}}
\def\pintegers{{\mathbb Z}_+}
\def\nintegers{{\mathbb N}}

\def\prb#1{\def\prbone{#1}
  \ifx\prbone\empty{\mathrm{P}}\else{\mathrm{P[\;}}#1{\mathrm{\;]}}\fi}
\def\prbseq#1#2{\def\prbseqone{#2}
  \ifx\prbseqone\empty{\mathrm{P}}_{#1}\ignorespaces
  \else{\mathrm{P}}_{#1}{\mathrm{[\;}}#2{\mathrm{\;]}}\fi}

\def\EEseq#1#2{\def\EEseqone{#2}
  \ifx\EEseqone\empty{\mathrm{E}}_{#1}\else
 {\mathrm{E}}_{#1}{\dsp\mathrm{[\;}}#2{\mathrm{\;]}}\fi}

\def\VVseq#1#2{\def\VVseqone{#2}
  \ifx\VVseqone\empty{\matrm{V}}_{#1}\else
 {\mathrm{V}}_{#1}{\dsp\mathrm{[\;}}#2{\mathrm{\;]}}\fi}

\def\img{\mathop{\mathrm{Im}}\nolimits}

\def\chrfcn#1{\mathop{\mathbf{1}}\nolimits_{#1}}

\def\sgcn{{\mathbb{\cal X}}}

\def\subsetne{\mbox{\raisebox{0.15em}{\setbox0=\hbox{$\subset$}\copy0\kern-0.5\wd0}\raisebox{-0.60em}{\setbox0=\hbox{$\ne$}\kern-0.5\wd0\copy0\kern-0.5\wd0}\setbox0=\hbox{$\subset$}\kern+0.5\wd0}}

\makeatletter
\def\address#1#2{\begingroup
\noindent\parbox[t]{7.8cm}{%
\small{\scshape\ignorespaces#1}\par\vskip1ex
\noindent\small{\itshape E-mail address}%
\/: #2\par\vskip4ex}\hfill%
\endgroup}%
\makeatother

\title{\uppercase{Choquet integration with submodular function on measurable space with sigma-algebra generating chain}} 
\author{
\bigskip \\
\textsc{
Tetsuya Hattori
} 
}
\date{} 
\begin{document}

\maketitle

\footnote{This work was supported by JSPS KAKENHI Grant Number 22K03358.}

\begin{abstract}
Based on a study of a formula representing submodular set function
as a supremum of measures dominated by the set function, 
we present a corresponding formula for a Choquet integration with respect to
the set function, on a measurable space which has a chain of measurable set
generating the sigma-algebra.
As an application we reproduce a basic formula in mathematical finance
on law invariant coherent risk measures.
We also study a recursion relation of set functions for which
the representation formula characterizes the fixed point.
\end{abstract}

\section{Introduction}
\seca{intro}

Let $(\Omega,{\cal F})$ be a measurable space,
namely, a $\sigma$-algebra ${\cal F}$ is a class of
subsets of $\Omega$ and is closed under 
complements and countable unions.
For a measurable set $A\in{\cal F}$ denote by ${\cal F}|_{A}$,
the class of measurable sets restricted to $A$, and
denote the set of finite measures on the measurable space $(A,{\cal F}|_A)$
by ${\cal M}(A)$.

For a real valued set function $v:\ {\cal F}\to\reals$ and
a measurable set $A\in{\cal F}$,
let ${\cal C}_{-,v}(A)$ be a class of measures dominated by $v$ on $A$;
\begin{equation}
\eqna{submodularlocalcoredef}
{\cal C}_{-,v}(A)=\{\mu\in{\cal M}(A)
\mid \mu(A)=v(A),
\ \mu(B)\le v(B),\ B\in {\cal F}|_{A}\}.
\end{equation}
If
\begin{equation}
\eqna{submodulartotalcorecoherentriskdef}
v(B)=\sup_{\mu\in{\cal C}_{-,v}(A)} \mu(B),
\end{equation}
holds for all $A,B\in{\cal F}$ satisfying $B\subset A$,
then it is easy to see that
\begin{equation}
\eqna{submodular}
v(A)+v(B)\ge v(A\cup B)+v(A\cap B), \ \ A,B\in {\cal F},
\end{equation}
holds \cite[Proposition 1]{cmod23}.
A set function which satisfies \eqnu{submodular} is called
a submodular function.

In \cite[Theorem 3]{cmod23} we gave a proof of a converse
that \eqnu{submodular} implies \eqnu{submodulartotalcorecoherentriskdef},
for a non-decreasing and continuous set function $v$,
when ${\cal F}$ is generated by a chain.
Here, by non-decreasing we mean
\begin{equation}
\eqna{nondecreasing}
v(A)\le v(B),\ A\subset B,\ A,B\in{\cal F}.
\end{equation}
We will later 	define the continuity of set functions which we adopt
in this paper,
and first make the assumption on ${\cal F}$ precise.
We consider, as in \cite{cmod23}, a following set of conditions
for a class of measurable sets ${\cal I}\subset{\cal F}$;
\begin{equation}
\eqna{sigmaalgebrawithtotallyorderedgeneratingclass}
\left\{\arrb{l}\dsp i)\ \ 
\sigma[{\cal I}]={\cal F},
\mbox{ where $\sigma[{\cal I}]$ denotes the smallest $\sigma$-algebra
 containing ${\cal I}$,}
\\ \dsp ii)\ \ 
\emptyset\in{\cal I},\ \Omega\in{\cal I},\\ \dsp iii)\ \ 
{\cal I}\mbox{ is a chain, i.e., totally ordered with respect to inclusion,}
\\ \dsp \phantom{iii) } \mbox{i.e., for all $I_1,I_2\in{\cal I}$
either $I_1\subset I_2$ or $I_2\subset I_1$\,.}
\arre\right.
\end{equation}
Let $\sgcn$ denote the collection of such classes;
\begin{equation}
\eqna{sgcn}
\arrb{l}\dsp 
\sgcn=\{{\cal I}\subset {\cal F}\mid
 \mbox{\eqnu{sigmaalgebrawithtotallyorderedgeneratingclass} holds} \}
\\ \dsp \phantom{\sgcn} =\{{\cal I}\subset {\cal F}\mid
\ {\cal I} \mbox{ is a chain such that }\ 
\emptyset\in{\cal I},\ \Omega\in{\cal I},\ \sigma[{\cal I}]={\cal F} \},
\arre
\end{equation}
and for ${\cal I}\in\sgcn$ and $A\in{\cal F}$,
define ${\cal I}_A$, the insertion of $A$ into ${\cal I}$, by
\begin{equation}
\eqna{insertion0}
{\cal I}_{A}
=\{ A\cap I \mid I\in{\cal I}\}\cup\{ A\cup I \mid I\in{\cal I}\}.
\end{equation}
Note that $\sgcn$ is closed under the insertion;
if ${\cal I}\in\sgcn$ and $A\in{\cal F}$ then 
$\dsp A\in {\cal I}_{A}\in\sgcn$ \cite[Lemma 2]{cmod23}.
With these notations, \cite[Proposition~1]{cmod23} and \cite[Theorem~3]{cmod23}
imply that if $\sgcn\ne\emptyset$, then a non-decreasing continuous 
set function $v$ is submodular, if and only if
\eqnu{submodulartotalcorecoherentriskdef} holds
for all $A,B\in{\cal F}$ satisfying $B\subset A$.
See \cite{DelbaenOsaka2008,Denneberg,71shapl} and 
the references in \cite{cmod23} for a part of 
large background and long history of the basic theory of
submodular functions.

It is proved in \cite[Proposition 7]{cmod23} that 
Polish spaces (separable, complete, metric space) with Borel $\sigma$-algebra
(the smallest $\sigma$-algebra containing all open balls) are
examples of measurable spaces satisfying $\sgcn\ne\emptyset$,
so that the above mentioned equivalence between submodular property
\eqnu{submodular} and the representation formula
\eqnu{submodulartotalcorecoherentriskdef} hold for Polish spaces.
One dimensional Borel $\sigma$-algebra and a finite set 
$\dsp \Omega_m=\{1,2,\ldots,m\}$ with $m$ elements and 
$\dsp {\cal F}=2^{\Omega_m}$ are 
among examples of Polish spaces, as well as a wide class of spaces 
extensively used in the theory of stochastic processes.
In this paper we exclusively consider $(\Omega,{\cal F})$
satisfying $\sgcn\ne\emptyset$,
and advance our study in \cite{cmod23} 
in the following directions.
\itdb
\item
The definition \eqnu{submodularlocalcoredef} of ${\cal C}_{-,v}(A)$
implies that
\eqnu{submodulartotalcorecoherentriskdef} follows from \eqnu{submodular} if 
there exists a measure $\mu\in {\cal C}_{-,v}(A)$ such that
$v(B)=\mu(B)$ holds.
Such a measure is studied by Shapley \cite{71shapl} as an extremal measure
when $\Omega$ is a finite set.
In our formulation, the extremal measure of $v$ is defined for each
${\cal I}\in \sgcn$.
(In \cite[Theorem~3]{cmod23}, for the sake of conciseness of the statement 
of the theorem, the notion of extremal measures is not explicitly written,
and is implicitly introduced in the proof of the theorem.)
For a non-decreasing set function $v$,
we define $\dsp\mu_{v,{\cal I}}:\ {\cal I}\to\preals$ by
\begin{equation}
\eqna{extremept}
 \mu_{v,{\cal I}}(I)=v(I),\ I\in{\cal I}.
\end{equation}
To extend $\dsp \mu_{v,{\cal I}}$ to a measure,
we define in \dfnu{vconti} in \secu{submodformula},
that a set function $v$ is continuous, if
for every ${\cal I}\in\sgcn$
$\dsp\mu_{v,{\cal I}}$\,, defined by \eqnu{extremept} on ${\cal I}$ 
and uniquely extended to the finite algebra ${\cal J}$
generated by ${\cal I}$ as a finitely additive measure, is
continuous on ${\cal J}$ in the standard measure theory sense.
(We will see in \secu{submodformula} that this definition is 
a generalization of the definition of continuity given in \cite{cmod23}
for submodular functions.)
If $v$ is continuous then $\dsp \mu_{v,{\cal I}}$ is uniquely extended to
a measure on ${\cal F}$, which we call the extremal measure of $v$ 
corresponding to ${\cal I}\in \sgcn$.
We prove in \secu{submodformula}
another representation formula for submodular function
in terms of the extremal measures
\begin{equation}
\eqna{submodextremeptrep}
 v=\sup_{{\cal I}\in \sgcn} \mu_{v,{\cal I}}
\end{equation}
which generalizes a corresponding result of theory of cores
for a finite set in \cite{71shapl} 
to measurable spaces satisfying $\sgcn\ne\emptyset$.

\item
We can define Choquet integration $v(f)$
with respect to non-decreasing continuous
set function $v$ for an integrable function $f$. 
If, in addition, $v$ is submodular,
the functional $\rho$ defined by $\rho(f)=v(-f)$ satisfies
the definition of the coherent risk measure studied in the field of
mathematical finance.
We then obtain a representation formula for $\rho$,
based on \eqnu{submodextremeptrep},
which reproduces a basic formula for a coherent risk measure with Fatou
property studied in mathematical finance \cite{risk99,risk01,riskmSK}.
See \cite{19Peng} for generalization and further development on 
Choquet integration with respect to submodular functions including
its relation to coherent risk measure.
In \cite{cmod23} we briefly announced an outline of 
how these results are formulated in our framework,
for which we give precise statements and proofs in \secu{Choquetint}.
As a further example in \secu{lawinvcrm}, we reproduce in our framework
a formula studied in \cite{riskmSK} for law-invariant coherent risk measures.

\item
Suggested by the representation formula 
\eqnu{submodextremeptrep} for submodular functions,
we consider the recursion relation 
\begin{equation}
\eqna{submodextremeptreprecursion}
 v_{n+1}=\sup_{{\cal I}\in \sgcn} \mu_{v_n,{\cal I}}\,,\ n=0,1,2,\ldots,
\end{equation}
with $v_0:\ {\cal F}\to\reals$ being a non-decrasing continuous set function.
Note that \eqnu{submodextremeptrep} 
implies $v_1\ne v_0$ if $v_0$ is not submodular and
$v_1=v_0$ if $v_0$ is submodular.
(We prove this equivalence in 
\thmu{totorder2modularvarrepbigcoreequiv} in \secu{submodformula}.)

For a finite set $(\Omega_m,2^{\Omega_m})$ all the non-decreasing
set functions are
continuous in our definition. We show in \secu{submodextremeptreprecursion}
that if $m=3$ and $v_0$ is non-decreasing then $v_1$ is submodular,
hence $v_1=v_2=\cdots$ in \eqnu{submodextremeptreprecursion},
while we show an example of $v_0$ for $m=4$ such that $v_1\ne v_2$\,.
\itde

Note that while we define and state results in terms of submodular functions,
the corresponding results in this paper for supermodular functions 
(convex games)
hold through a well-known correspondence
\begin{equation}
\eqna{convexconcavedualtr}
\tilde{v}(A)=v(\Omega)-v(A^c)+v(\emptyset),\ A\in{\cal F},
\end{equation}
which gives a non-decreasing continuous supermodular (resp., submodular)
function $\tilde{v}$ from a  non-decreasing continuous submodular
(resp., supermodular) function $v$
satisfying $\tilde{v}(\Omega)=v(\Omega)$ and 
$\tilde{v}(\emptyset)=v(\emptyset)$.
We will choose definitions and statements
in \secu{submodformula}
so that corresponding results for convex games also hold by 
the correspondence \eqnu{convexconcavedualtr}.
In fact, classical theories by Shapley \cite{71shapl}, as well as
many works on cooperative game theories, are written in terms of
convex games, for which we should use \eqnu{convexconcavedualtr}
when comparing with the results for sub-modular functions.

For notational simplicity we assume $v(\emptyset)=0$ for any set function $v$
throughout this paper. The formula in this paper can be generalized to the case
$v(\emptyset)\ne 0$ by the replacements $v(A)\mapsto v(A)-v(\emptyset)$.

\paragraph*{Acknowledgement.}
The author deeply thanks the referee for the valuable comments,
including (but not restricted to) the proof of 
submodularity and \eqnu{upperconti} for the limit in \prpu{recursionconv},
and bringing the reference \cite{19Peng} to the author's attention.

\section{Representation formula for submodular function}
\seca{submodformula}

In this section, we fix a measurable space $(\Omega,{\cal F})$ and
assume $\sgcn\ne\emptyset$, where $\sgcn$ is defined in \eqnu{sgcn}
in \secu{intro}.

Let $v:\ {\cal F}\to\reals$ be a non-decreasing set function.
Let ${\cal I}\in\sgcn$, and denote by ${\cal J}$
the finitely additive class generated by $\dsp {\cal I}$\,,
the smallest class of sets such that ${\cal I}$ is a subset and
closed under complement and union.
Since $\dsp {\cal I}$ is totally ordered with respect to inclusion,
we have an explicit representation
\begin{equation}
\eqna{finiteadditiveclass}
\arrb{l}\dsp
{\cal J}=
\{ \bigcup_{i=1}^n (C_i\cap D_i^c) \mid
C_1\supset D_1\supset C_2\supset\cdots\supset D_n\,,
\\ \dsp \phantom{{\cal J}=\{ }
\ C_i,\ D_i\in{\cal I}\,,\ i=1,2,\ldots,n,\ \ n=1,2,3,\ldots \}.
\arre
\end{equation}
Using the notation in the right-hand side of \eqnu{finiteadditiveclass},
we can define, as in an elementary textbook on measure theory
(see \cite{cmod23}),
a finitely additive set function $\dsp\mu_{v,{\cal I}}:\ {\cal J}\to\reals$ 
as a unique extension of \eqnu{extremept}, by
\begin{equation}
\eqna{finiteadditivemeasure}
\mu_{v,{\cal I}}(\bigcup_{i=1}^n (C_i\cap D_i^c))=
\sum_{i=1}^n (v(C_i)-v(D_i)).
\end{equation}

%
\dfnb
\dfna{vconti}
We say that a non-decreasing set function
$v:\ {\cal F}\to\reals$ is continuous, if 
for every ${\cal I}\in\sgcn$ the finitely additive measure
$\dsp \mu_{v,{\cal I}}:\ {\cal J}\to\reals$ defined by
\eqnu{finiteadditivemeasure} is $\sigma$-additive, or equivalently,
continuous, on the finitely additive class ${\cal J}$.

A standard extension theorem of measures and $\sigma[{\cal I}]={\cal F}$
in \eqnu{sgcn} then imply that $\dsp \mu_{v,{\cal I}}$ is uniquely
extended to a measure on ${\cal F}$ for all ${\cal I}\in\sgcn$ 
if $v$ is continuous.
We will use the same notation and call the measure
$\dsp \mu_{v,{\cal I}}:\ {\cal F}\to\reals$ the 
extremal measure of $v$ corresponding to ${\cal I}\in \sgcn$.
\DDD
\dfne
We note that this definition of continuity of $v$ is consistent
with the definition of continuity in \cite{cmod23} for 
submodular functions. Namely, the following holds.
\prpb
\prpa{conti24equivconti23onsubmod}
Let $(\Omega,{\cal F})$ be a measurable space satisfying $\sgcn\ne\emptyset$
and $v:\ {\cal F}\to\reals$ be a non-decreasing submodular function
satisfying $v(\emptyset)=0$,
i.e., a set function satisfying $v(\emptyset)=0$, 
\eqnu{nondecreasing} and \eqnu{submodular}.
Then $v$ is continuous if and only if
\begin{equation}
\eqna{upperconti}
\limf{n} v(A_n)=v(\bigcup_{n\in\nintegers} A_n),
\ \ A_1\subset A_2\subset\cdots,\ \ A_n\in{\cal F},\ n=1,2,3,\ldots,
\end{equation}
and
\begin{equation}
\eqna{lowerconti}
\limf{n} v(A_n)=v(\bigcap_{n\in\nintegers} A_n)
\ \ A_1\supset A_2\supset\cdots,\ \ A_n\in{\cal F},\ n=1,2,3,\ldots,
\end{equation}
hold.
\DDD
\prpe
We note that our definition of continuity implies
\eqnu{upperconti} and \eqnu{lowerconti} without an assmption of
submodularity.
\prpb
\prpa{conti24toconti23}
Let $(\Omega,{\cal F})$ be a measurable space satisfying $\sgcn\ne\emptyset$
and $v:\ {\cal F}\to\reals$ be a non-decreasing continuous set function
satisfying $v(\emptyset)=0$.
Then \eqnu{upperconti} and \eqnu{lowerconti} hold.
\DDD
\prpe
To prove \prpu{conti24equivconti23onsubmod} and \prpu{conti24toconti23}
(as well as for later use) we prepare the following 
\lemu{terminalmeasureinC} and \lemu{ctblechaininsertion}.
\lemb
\lema{terminalmeasureinC}
Let $(\Omega,{\cal F})$ be a measurable space satisfying $\sgcn\ne\emptyset$,
and let ${\cal I}\in\sgcn$. Denote by ${\cal J}$ the finitely additive class
\eqnu{finiteadditiveclass} generated by ${\cal I}$.

If a non-decreasing set function $v:\ {\cal F}\to\reals$ satisfies
\eqnu{submodular}, i.e., submodular, then the finitely additive measure
$\dsp \mu_{v,{\cal I}}:\ {\cal J}\to\reals$ definedd by
\eqnu{finiteadditivemeasure} satisfies
$\dsp \mu_{v,{\cal I}}(J)\le v(J)$, $J\in {\cal J}$.
\DDD
\leme
\prfb
A proof is basically same as that of \cite[Lemma~6]{cmod23}.
Using the expression \eqnu{finiteadditiveclass} let
\[
\arrb{l}\dsp
J=\bigcup_{i=1}^{n} (C_i\cap (D_i)^c)\in {\cal J};\ \ 
C_1\supset D_1\supset C_2\supset\cdots\supset D_{n}\,,
\\ \dsp
\ C_i,\ D_i\in{\cal I}\,,\ i=1,2,\ldots,n,
\arre
\]
and put
$\dsp
A_i=\bigcup_{j=i}^{n} (C_j\cap (D_{j})^c)$, $i=1,2,\ldots,n$,
and $A_{n+1}=\emptyset$.
Then $A_1=J$ and
\[
A_i\cup D_i=C_{i}\,,\ A_i\cap D_i=A_{i+1}\,,\ \ i=1,2,\ldots,n,
\]
Apply \eqnu{submodular} with $A=A_i$ and $B=D_i$ to find
\[
 v(C_{i})+v(A_{i+1})\le v(A_i)+v(D_i),\ i=1,2,\ldots,n.
\]
Summing this up with $i$ and using
\eqnu{finiteadditivemeasure} leads to
$\dsp\mu_{v,{\cal I}}(J)\le v(J)$,
\QED
\prfe
Let $(\Omega,{\cal F})$ be a measurable space satisfying $\sgcn\ne\emptyset$,
and let ${\cal I}\in\sgcn$. 

Before moving on to the next Lemma, we remark on sequential insertion to
${\cal I}$.
For $A,B\in{\cal F}$, a sequential insertion of $A$ and $B$ to ${\cal I}$
in general depends on the order of the insertion. In fact,
\begin{equation}
\eqna{insertionAB}
\arrb{l}\dsp 
({\cal I}_A)_B
=\{I'\cap B\mid I'\in{\cal I}_A\}\cup\{I'\cup B\mid I'\in{\cal I}_A\}
\\ \dsp \phantom{({\cal I}_A)_B}
=\{I\cap (A\cap B)\mid I\in{\cal I}\}
\cup \{(I\cap B)\cup (A\cap B)\mid I\in{\cal I}\}
\\ \dsp \phantom{({\cal I}_A)_B=}
\cup \{(I\cup B)\cap (A\cup B)\mid I\in{\cal I}\}
\cup \{I\cup (A\cup B)\mid I\in{\cal I}\}
\arre
\end{equation}
implies $\dsp B\in ({\cal I}_A)_B$\,,
while $A$ may not be an element.
Note, however, that \eqnu{insertionAB} implies
$\dsp A\cap B,\ A\cup B\in ({\cal I}_A)_B \cap ({\cal I}_B)_A$\,.

If $\{A,B\}$ is a chain, i.e., $A\subset B$ or $A\supset B$, 
then the order of insertion is irrelevant, and we will use the notation
$\dsp {\cal I}_{A,B}:=({\cal I}_A)_B=({\cal I}_B)_A$\,.
For example, if $B\subset A$ then
\begin{equation}
\eqna{chaininsertion}
{\cal I}_{A,B}:=({\cal I}_A)_B=\{I\cap B,\ (I\cup B)\cap A,\ I\cup A\mid
I\in{\cal I}\}=({\cal I}_B)_A\in\sgcn.
\end{equation}

For ${\cal A}=\{A_n\mid n=1,2,3,\ldots\}\subset{\cal F}$ satisfying
$A_1\subset A_2\subset\cdots$, define the (countable) insertion
$\dsp {\cal I}_{\cal A}$ of ${\cal A}$ into ${\cal I}$ by
\begin{equation}
\eqna{ctblechaininsertion1}
{\cal I}_{\cal A}=
 \{I \cap A_1\mid I\in{\cal I}\}
\cup\{(I \cup A_n) \cap A_{n+1}\mid I\in{\cal I},\ n\in\nintegers\}
\cup\{ I\cup \bigcup_{n\in\nintegers} A_n\mid I\in{\cal I}\}.
\end{equation}
Similarly,
For ${\cal A}=\{A_n\mid n=1,2,3,\ldots\}\subset{\cal F}$ satisfying
$A_1\supset A_2\supset\cdots$, define
\begin{equation}
\eqna{ctblechaininsertion2}
{\cal I}_{\cal A}=
 \{ I\cap \bigcap_{n\in\nintegers} A_n\mid I\in{\cal I}\}
\cup\{(I \cup A_{n+1}) \cap A_{n}\mid I\in{\cal I},\ n\in\nintegers\}
\cup \{I \cup A_1\mid I\in{\cal I}\}.
\end{equation}
\lemb
\lema{ctblechaininsertion}
$\dsp {\cal I}_{\cal A}$ of \eqnu{ctblechaininsertion1} satisfies
$\dsp A_n\in{\cal I}_{\cal A}$\,, $n\in\nintegers$,
$\dsp \bigcup_{n\in\nintegers} A_n\in{\cal I}_{\cal A}$\,, and 
$\dsp {\cal I}_{\cal A}\in\sgcn$.

$\dsp {\cal I}_{\cal A}$ of \eqnu{ctblechaininsertion2} satisfies
$\dsp A_n\in{\cal I}_{\cal A}$\,, $n\in\nintegers$,
$\dsp \bigcap_{n\in\nintegers} A_n\in{\cal I}_{\cal A}$\,, and 
$\dsp {\cal I}_{\cal A}\in\sgcn$.
\DDD
\leme
\prfb
Consider first the case $A_1\subset A_2\subset\cdots$.
Since ${\cal I}\in\sgcn$ implies that ${\cal I}$ is a chain containing
$\emptyset$ and $\Omega$, 
$\dsp {\cal I}_{\cal A}$ also shares these properties, because,
for example, if $m<n$ then for $I,I'\in{\cal I}$,
\[
(I\cup  A_m)\cap A_{m+1}\subset A_{m+1}\subset A_n\subset
(I'\cup  A_n)\cap A_{n+1}.
\]
To prove $\dsp{\cal I}_{\cal A}\in\sgcn$, it only remains to prove
$\dsp \sigma[{\cal I}_{\cal A}]={\cal F}$.

Before proceeding with proving this, note that 
$\dsp \bigcup_{n\in\nintegers} A_n=\emptyset\cup \bigcup_{n\in\nintegers} A_n
\in{\cal I}_{\cal A}$,
$\dsp A_1=\Omega\cap A_1\in{\cal I}_{\cal A}$,
and for $n\in\nintegers$,
$\dsp A_{n+1}=(\Omega\cup A_n)\cap A_{n+1}\in{\cal I}_{\cal A}$,

Returning to the proof of 
$\dsp \sigma[{\cal I}_{\cal A}]={\cal F}$,
note that ${\cal I}\in\sgcn$ implies $\sigma[{\cal I}]={\cal F}$, hence
it suffices to prove $\sigma[{\cal I}_{\cal A}]\supset {\cal I}$.
To prove this, let $I\in{\cal I}$, and note that
\[
I=(I\cap (\bigcup_{n\in\nintegers} A_n)^c)\cup
(I\cap A_1)\cup \bigcup_{n\in\nintegers}(I\cap A_n^c\cap A_{n+1}).
\]
Since
$\dsp I\cup \bigcup_{n\in\nintegers} A_n\in {\cal I}_{\cal A}$ and
$\dsp \bigcup_{n\in\nintegers} A_n\in{\cal I}_{\cal A}$\,,
\[
I\cap (\bigcup_{n\in\nintegers} A_n)^c)=
(I\cup \bigcup_{n\in\nintegers} A_n)\cap (\bigcup_{n\in\nintegers} A_n)^c
\in \sigma[{\cal I}_{\cal A}]\,,
\]
and for $n\in\nintegers$ since
$\dsp (I \cup A_n) \cap A_{n+1}\in{\cal I}_{\cal A}$ and
$\dsp A_n\in {\cal I}_{\cal A}$\,,
\[
I\cap A_n^c\cap A_{n+1}=
((I \cup A_n) \cap A_{n+1})\cap A_n^c\in\sigma[{\cal I}_{\cal A}].
\]
Therefore $\dsp I\in \sigma[{\cal I}_{\cal A}]$,
which proves $\dsp {\cal I}\subset \sigma[{\cal I}_{\cal A}]$.

Next let $A_1\supset A_2\supset\cdots$.
All the claims except $\sigma[{\cal I}_{\cal A}]\supset {\cal I}$
are proved in a similar way as the previous case.
Note that for $I\in{\cal I}$ we have
\[
I=(I\cap \bigcap_{n\in\nintegers} A_n)\cup
(I\cap A_1^c)\cup \bigcup_{n\in\nintegers}(I\cap A_{n+1}^c\cap A_{n}).
\]
Since $I\cup A_1\in{\cal I}_{\cal A}$ and
$A_1=\emptyset\cup A_1\in{\cal I}_{\cal A}$,
\[
I\cap A_1^c=(I\cup A_1)\cap A_1^c\in\sigma[{\cal I}_{\cal A}],
\]
and for $n\in\nintegers$ since
$\dsp (I \cup A_{n+1}) \cap A_{n}\in{\cal I}_{\cal A}$ and
$\dsp A_{n+1}\in {\cal I}_{\cal A}$\,,
\[
I\cap A_{n+1}^c\cap A_{n}=
((I \cup A_{n+1}) \cap A_{n})\cap A_{n+1}^c\in\sigma[{\cal I}_{\cal A}].
\]
Therefore $\dsp I\in \sigma[{\cal I}_{\cal A}]$,
which proves $\dsp {\cal I}\subset \sigma[{\cal I}_{\cal A}]$.
\QED
\prfe
\prfofb{\protect\prpu{conti24toconti23}}
Let $v$ be a non-decreasing continuous set function
satisfying $v(\emptyset)=0$.
To prove \eqnu{lowerconti}, let
$\dsp{\cal A}=\{A_n\mid n=1,2,3,\ldots\}\subset{\cal F}$ 
be a sequence satisfying $A_1\supset A_2\supset\cdots$,
and let $\dsp{\cal I}_{\cal A}$ be as in \eqnu{ctblechaininsertion2}.
\lemu{ctblechaininsertion} then implies 
$\dsp{\cal I}_{\cal A}\in\sgcn$.
Denote by $\dsp{\cal J}_{\cal A}$ the finitely additive class
generated by $\dsp{\cal I}_{\cal A}$.
Since $v$ is continuous, \dfnu{vconti} implies that
$\dsp \mu_{v,{\cal I}_{\cal A}}:\ {\cal J}_{\cal A}\to\reals$ is
continuous.
Since \lemu{ctblechaininsertion} implies 
$\dsp A_n\in{\cal I}_{\cal A}\subset{\cal J}_{\cal A}$, $n\in\nintegers$,
and $\dsp \bigcap_{n\in\nintegers} A_n
\in{\cal I}_{\cal A}\subset{\cal J}_{\cal A}$,
continuity of $\dsp  \mu_{v,{\cal I}_{\cal A}}$ on
$\dsp{\cal J}_{\cal A}$ implies
$\limf{n} \mu_{v,{\cal I}_{\cal A}}(A_n)
=\mu_{v,{\cal I}_{\cal A}}(\bigcap_{n\in\nintegers} A_n)$.
Also \eqnu{extremept} implies
$\dsp \mu_{v,{\cal I}_{\cal A}}(A_n)=v(A_n)$, $n\in\nintegers$,
and $\dsp \mu_{v,{\cal I}_{\cal A}}(\bigcap_{n\in\nintegers} A_n)
=v(\bigcap_{n\in\nintegers} A_n)$.
Therefore we have
$\limf{n} v(A_n)=v(\bigcap_{n\in\nintegers} A_n)$,
which proves \eqnu{lowerconti}.

Proof of \eqnu{upperconti} is similar to that of \eqnu{lowerconti},
if we replace $A_1\supset A_2\supset\cdots$,
$\dsp \bigcap_{n\in\nintegers} A_n$\,,
and \eqnu{ctblechaininsertion2}
by $A_1\subset A_2\subset\cdots$,
$\dsp \bigcup_{n\in\nintegers} A_n$\,,
and \eqnu{ctblechaininsertion1}, respectively.
\QED
\prfofe
\prfofb{\protect\prpu{conti24equivconti23onsubmod}}
That the continuity of $v$ implies \eqnu{upperconti} and \eqnu{lowerconti}
is a consequence of \prpu{conti24toconti23}.
To prove the converse, assume that
$v$ satisfies \eqnu{submodular}, \eqnu{nondecreasing},
\eqnu{upperconti}, and \eqnu{lowerconti}.
Let ${\cal I}\in\sgcn$ and ${\cal J}$ be the finitely additive class
generated by ${\cal I}$.
As in standard measure theory, to prove continuity of
$\dsp \mu_{v,{\cal I}}$ on ${\cal J}$ it suffices to
prove that for any sequence $A_n\in{\cal J}$, $n=1,2,3,\ldots$,
satisfying $A_1\supset A_2\supset \cdots$ and 
$\dsp \bigcap_{n\in\nintegers} A_n=\emptyset$,
$\dsp \limf{n}\mu_{v,{\cal I}}(A_n)=0$.
For such sequence $\{A_n\}$, \eqnu{lowerconti} implies,
with $v(\emptyset)=0$,
$\dsp\limf{n}v(A_n)=0$.
On the other hand \lemu{terminalmeasureinC} implies 
$\dsp\mu_{v,{\cal I}}(A_n)\le v(A_n)$, $n\in\nintegers$.
These with non-negativity of measure $\dsp\mu_{v,{\cal I}}$
imply $\dsp \limf{n}\mu_{v,{\cal I}}(A_n)=0$,
which completes the proof of continuity of $v$.
\QED
\prfofe
\remb
In \prpu{conti24equivconti23onsubmod} both 
\eqnu{upperconti} and \eqnu{lowerconti} are stated in order to make
the statement also hold when `submodular' is replaced by `supermodular'.
It is a known elementary fact that if $v$ is non-decreasing and submodular 
then \eqnu{lowerconti} implies \eqnu{upperconti}.

Note that if $v$ is non-decreasing submodular and \eqnu{upperconti} holds,
\eqnu{lowerconti} does not necessarily follow.
In fact, if $\dsp \mu:\ {\cal F}\to\reals$ is a finite measure such that
there exists a sequence $\dsp A_n\in{\cal F}$, $n\in\nintegers$,
such that 
$\dsp A_1\supset A_2\supset\cdots$, $\dsp \mu(A_n)>0$, $n\in\nintegers$,
and $\dsp \bigcap_{n\in\nintegers} A_n=\emptyset$,
then $\dsp v:\ {\cal F}\to\reals$ defined by
$\dsp v(A)=\left\{\arrb{ll} \mu(A)+1, & A\ne \emptyset,\\ 0 & A=\emptyset,
\arre\right.$
is non-decreasing submodular and \eqnu{upperconti} holds,
but \eqnu{lowerconti} fails because
\[
\limf{n} v(A_n)=\limf{n} \mu(A_n)+1= \mu(\bigcap_{n\in\nintegers} A_n)+1
=1\ne 0=v(\bigcap_{n\in\nintegers} A_n).
\]
An explicit example is the one-dimensional Lebesgue measure $\mu=\mu_1$
on an interval $(\Omega,{\cal F},\mu)=((0,1],{\cal B}_1((0,1]),\mu_1)$
and $\dsp A_n=(0,\frac1n]$, $n\in\nintegers$.
\DDD\reme

We move on to the main theorem of this section on the
representation formula of submodular functions.
\thmb
\thma{totorder2modularvarrepbigcoreequiv}
Let $(\Omega,{\cal F})$ be a measurable space satisfying $\sgcn\ne\emptyset$.
Let $v:\ {\cal F}\to\reals$ be a non-decreasing and continuous 
set function, satisfying $v(\emptyset)=0$,
where continuity of set function is as in \dfnu{vconti}.
Then the following (a), (b), (c), (d) are equivalent.
\def\labelenumi{(\alph{enumi})}
\itmb
\item $v$ is submodular, i.e., \eqnu{submodular} holds.

\item
For all ${\cal I}\in\sgcn$, the extremal measure 
$\dsp \mu_{v,{\cal I}}$ of $v$ corresponding to ${\cal I}$
determined by \eqnu{extremept} satisfies
$\dsp \mu_{v,{\cal I}}\in {\cal C}_{-,v}(\Omega)$,
where $\dsp {\cal C}_{-,v}(\Omega)$ is defined in
\eqnu{submodularlocalcoredef} with $A=\Omega$.

\item
The representation \eqnu{submodextremeptrep} holds.

\item
For all $A,B\in{\cal F}$ satisfying $B\subset A$,
\eqnu{submodulartotalcorecoherentriskdef} holds.
\DDD
\itme
\def\labelenumi{(\roman{enumi})}
\thme
\prfb
We will prove
(a) $\Rightarrow$ (b) $\Rightarrow$ (c) $\Rightarrow$ (d)  $\Rightarrow$ (a) 
in this order.

\paragraph{(a) $\Rightarrow$ (b).}
Let ${\cal I}\in\sgcn$ and let ${\cal J}$ denote the finitely additive class 
generated by ${\cal I}$.
Since $\Omega\in{\cal I}$, $\dsp \mu_{v,{\cal I}}(\Omega)= v(\Omega)$.
Assumption (a) and \lemu{terminalmeasureinC} imply
$\dsp \mu_{v,{\cal I}}(J)\le v(J)$, $J\in {\cal J}$.
By the assumption that $v$ is continuous, $\dsp \mu_{v,{\cal I}}$ is
uniquely extended to a measure on ${\cal F}$, which we also write
$\dsp \mu_{v,{\cal I}}$.
According to standard measure theory, the measure $\dsp \mu_{v,{\cal I}}$
on ${\cal F}=\sigma[{\cal I}]=\sigma[{\cal J}]$ is approximated by
the restriction of $\dsp \mu_{v,{\cal I}}$ to ${\cal J}$.
Hence $\dsp \mu_{v,{\cal I}}(A)\le v(A)$ holds for all $A\in {\cal F}$,
which further implies $\dsp\mu_{v,{\cal I}}\in{\cal C}_{-,v}(\Omega)$.

\paragraph{(b) $\Rightarrow$ (c)}.
Fix $A\in{\cal F}$ arbitrarily.
Since by assumption $\sgcn\ne\emptyset$, there is ${\cal I}\in \sgcn$.
Then the insertion of $A$ also satisfies $\dsp {\cal I}_A\in\sgcn$.
With $\dsp A\in {\cal I}_A$ we also have $\dsp\mu_{v,{\cal I}_A}(A)= v(A)$.
Therefore we have $\dsp v(A)\le \sup_{{\cal I}\in \sgcn} \mu_{v,{\cal I}}(A)$.
On the other hand, by assumption (b) and the definition of 
$\dsp {\cal C}_{-,v}(\Omega)$, we have
$\dsp\mu_{v,{\cal I}}(A)\le v(A)$ for all 
${\cal I}\in\sgcn$, which further implies
$\dsp \sup_{{\cal I}\in \sgcn} \mu_{v,{\cal I}}(A)\le v(A)$.
Hence the equality holds. Since $A\in{\cal F}$ is arbitrary,
\eqnu{submodextremeptrep} follows.

\paragraph{(c) $\Rightarrow$ (d)}.
Let $A,B\in{\cal F}$ satisfy $B\subset A$.
The definition of $\dsp {\cal C}_{-,v}(A)$ implies
$\mu(B)\le v(B)$ for all $\dsp\mu\in{\cal C}_{-,v}(A)$.
Hence to prove \eqnu{submodulartotalcorecoherentriskdef}
it suffices to prove existence of $\dsp\mu\in{\cal C}_{-,v}(A)$
satisfying $\dsp v(B)=\mu(B)$.

Note that $B\subset A$ and let $\dsp {\cal I}_{A,B}$ the 
sequential insertion of $A$ and $B$ in \eqnu{chaininsertion}.
Assumption (c) then implies $\dsp \mu_{v,{\cal I}_{A,B}}(C)\le v(C)$
for all $C\in{\cal F}$ satisfying $C\subset A$.
$\dsp A,B\in{\cal I}_{A,B}$ and \eqnu{extremept} imply
$\dsp\mu_{v,{\cal I}_{A,B}}(A)=v(A)$ and $\dsp\mu_{v,{\cal I}_{A,B}}(B)=v(B)$.
Hence the restriction $\dsp\mu_{v,{\cal I}_{A,B}}|_A$ of
$\dsp\mu_{v,{\cal I}_{A,B}}$ to $A$ is in $\dsp {\cal C}_{-,v}(A)$,
which is the measure we are looking for.

\paragraph{(d) $\Rightarrow$ (a)}
This claim is essentially \cite[Proposition 1]{cmod23}, whose proof
works here.
Let $A,B\in {\cal F}$ and 
substitute $A$ in \eqnu{submodulartotalcorecoherentriskdef} by
$A\cup B$ to obtain
$\dsp v(B)=\sup_{\mu\in{\cal C}_{-,v}(A\cup B)} \mu(B)$
and $\dsp v(A)=\sup_{\mu\in{\cal C}_{-,v}(A\cup B)} \mu(A)$,
which further imply $v(B)\ge \mu(B)$ and $v(A)\ge \mu(A)$ for 
all $\mu\in {\cal C}_{-,v}(A\cup B)$.
Also \eqnu{submodularlocalcoredef} implies
$v(A\cup B)=\mu(A\cup B)$ for all $\mu\in {\cal C}_{-,v}(A\cup B)$.
Finally, for $\eps>0$,
$\dsp v(A\cap B)=\sup_{\mu\in{\cal C}_{-,v}(A\cup B)} \mu(A\cap B)$
implies that there exists $\mu\in{\cal C}_{-,v}(A\cup B)$ such that
$\dsp v(A\cap B)\le \mu(A\cap B)+\eps$.
These equality and inequalities imply
\[\arrb{l}\dsp
v(A\cup B)+v(A\cap B)-v(A)-v(B) 
\\ \dsp {}
\le \mu(A\cup B)+\mu(A\cap B)+\eps-\mu(A)-\mu(B)
=\eps.
\arre\]
$\eps$ can be any positive constant, hence \eqnu{submodular} follows.
\QED
\prfe
%
%
%

In \cite{cmod23} we used the term continuous to mean 
\eqnu{upperconti} and \eqnu{lowerconti} 
also for non-decreasing supermodular functions (convex games). 
As noted at the end of \secu{intro},
\eqnu{convexconcavedualtr} implies that
\prpu{conti24equivconti23onsubmod} holds also for supermodular functions.
Furthermore, a convex game version of
\thmu{totorder2modularvarrepbigcoreequiv} also holds.
\corb
\cora{totorder2modularvarrepbigcoreequiv}
Let $(\Omega,{\cal F})$ be a measurable space satisfying $\sgcn\ne\emptyset$.
Let $v:\ {\cal F}\to\reals$ be a non-decreasing and continuous 
set function, satisfying $v(\emptyset)=0$.
Then the following (a), (b), (c), (d) are equivalent.
\def\labelenumi{(\alph{enumi})}
\itmb
\item $v$ is supermodular (convex game), namely,
$\dsp v(A)+v(B)\le v(A\cup B)+v(A\cap B), \ \ A,B\in {\cal F}$,
holds.

\item
For all ${\cal I}\in\sgcn$, the extremal measure 
$\dsp \mu_{v,{\cal I}}$ of $v$ corresponding to ${\cal I}$
determined by \eqnu{extremept} satisfies
$\dsp \mu_{v,{\cal I}}\in {\cal C}_{+,v}(\Omega)$,
where
$\dsp {\cal C}_{+,v}(\Omega)=\{\mu\in{\cal M}(\Omega)
\mid \mu(\Omega)=v(\Omega),\ \mu(B)\ge v(B),\ B\in {\cal F}\}$.

\item
It holds that
$\dsp v(A)=\inf_{{\cal I}\in \sgcn} \mu_{v,{\cal I}}(A),\ A\in{\cal F}$.

\item
For all $A,B\in{\cal F}$ satisfying $B\subset A$,
$\dsp v(B)=\inf_{\mu\in{\cal C}_{+,v}(A)} \mu(B)$ holds.
\DDD
\itme
\def\labelenumi{(\roman{enumi})}
\core

\section{Choquet integration and representation formula for submodular function}
\seca{Choquetint}
In this section, we fix a measurable space $(\Omega,{\cal F})$ which
satisfies $\sgcn\ne\emptyset$.

\subsection{Choquet integration}

Let $v:\ {\cal F}\to\reals$ be a non-decreasing set function,
satisfying $v(\emptyset)=0$, and $f:\ \Omega\to\reals$ a 
real valued measurable function.

Note that if we define $\dsp g:\ \reals\to\reals$ by 
$\dsp g(z)=v(\{\omega\in\Omega\mid f(\omega)>z\})$, 
and $\dsp g^{-1}(a)=\sup_{g(z)>a}z$ for $a\ge 0$,
non-decreasing property of $v$ implies 
$\dsp\{z\in\reals\mid g(z)>a\} = (-\infty,g^{-1}(a))$,
hence $g$ is a $1$-dimensional Borel measurable function,
and we can consider the standard $1$-dimensional
Lebesgue integration of $g$.
We define 
\begin{equation}
\eqna{Choquetgen}
v(f)=\lim_{y\to-\infty} \biggl(y\,v(\Omega)
+\int_{y}^{\infty} v(\{\omega\in\Omega\mid f(\omega)>z\})\,dz\biggr)
\end{equation}
whenever the right-hand side is a real value.
If either the Lebesgue integration or the limit diverges
in the right-hand side of \eqnu{Choquetgen} we do not define $v(f)$.
If $v(f)$ of \eqnu{Choquetgen} is defined it is equal to 
the asymmetric integral in terms of \cite[Chap.~5]{Denneberg}.
We will refer to \eqnu{Choquetgen} as Choquet integration 
(of $f$ with respect to $v$) in this paper.

If $f(\omega)\ge x$, $\omega\in\Omega$, holds for an $x\in \reals$,
then 
\[
\int_y^x v(\{\omega\in\Omega\mid f(\omega)>z\})dz
=(x-y)\,v(\Omega),\ \ y\le x,
\]
hence in this case \eqnu{Choquetgen} has a simpler expression
\begin{equation}
\eqna{Choquetbd}
v(f)=x\,v(\Omega)
+\int_{x}^{\infty} v(\{\omega\in\Omega\mid f(\omega)>z\})\,dz.
\end{equation}

Simple function, a measurable function whose image is a finite set,
is an elementary example of Choquet integrable function.
We find it convenient to adopt a representation
\begin{equation}
\eqna{simplefcnrepequiv2}
f=\sum_{i=1}^na_i\chrfcn{A_i}\,,
\end{equation}
where
\begin{equation}
\eqna{simplefcnChoquetrep}
\arrb{l}\dsp
a_i\ge0,\ i=1,\ldots,n-1,\ a_n\in\reals,
\\ \dsp \emptyset =A_0\subset 
A_1\subset A_2\subset\cdots\subset A_{n-1}\subset A_n=\Omega.
\arre
\end{equation}
Here and in the following, 
we use a symbol $\dsp\chrfcn{A}$ to denote a characteristic function
of a set $A$ 
defined by $\dsp \chrfcn{A}(\omega)=\left\{\arrb{ll} 1,& \omega\in A,\\
0,& \omega\not\in A.\arre\right.$
(We use $\dsp \chrfcn{\cdot}$ on any space.)

We then have
$\dsp A_i=\{\omega\in\Omega\mid f(\omega)\ge \sum_{j=i}^n a_j\}\in {\cal F},
\ i=1,\ldots,n,$
hence \eqnu{Choquetbd} with $x=a_n$ and \eqnu{simplefcnChoquetrep} implies
\begin{equation}
\eqna{simplefcnChoquetint}
v(f)=\sum_{i=1}^{n} a_i\,v(A_i).
\end{equation}
Note that $f$ of \eqnu{simplefcnrepequiv2} takes values 
$\dsp b_i=\sum_{j=i}^n a_j$\,, $i=1,\ldots,n$, with
$b_1\ge b_2\ge \cdots\ge b_n=a_n$\,.
With $\dsp B_i=A_i\cap A_{i-1}^c$, 
$i=1,\ldots,n$, we have
$B_i\cap B_j=\emptyset$, $i\ne j$, $\dsp \bigcup_{i=1}^n B_i=A_n=\Omega$, and
$\dsp f=\sum_{i=1}^nb_i\chrfcn{B_i}$, which perhaps is closer to
a familiar expression of a simple function.
With this expression, \eqnu{simplefcnChoquetint} implies
$\dsp v(f)=\sum_{i=1}^{n-1} (b_i-b_{i+1})\,v(\bigcup_{j=1}^iB_j)+b_nv(\Omega)$.

\subsection{Monotone convergene theorem}

Hereafter we assume that $v:\ {\cal F}\to\reals$ is 
a continuous and non-decreasing set function,
We then have an analog of monotone convergence theorem
for the Choquet integration with respect to $v$.
To state the theorem, 
we extend some elementary terms in measure theory to
continuous, non-decreasing set function $v$.
We say that a measurable set $N\in{\cal F}$ is a $v$-null set, if
\begin{equation}
\eqna{nullset}
v(N\cup A)=v(N^c\cap A)=v(A),\ A\in{\cal F}.
\end{equation}
For a sequence of measurable functions 
$f_n:\ \Omega\to\reals$, $n\in\nintegers$, and a 
measurable function $f:\ \Omega\to\reals$,
we say that the sequence converges $v$-almost surely to $f$, if
there exists a $v$-null set $N\in{\cal F}$ such that
\begin{equation}
\eqna{asconvergence}
\limf{n} f_n(\omega)=f(\omega),\ \omega\in N^c.
\end{equation}
\thmb
\thma{monconvthm}
Let $(\Omega,{\cal F})$ be a measurable space satisfying $\sgcn\ne\emptyset$,
and $v:\ {\cal F}\to\reals$ be a non-decreasing, continuous set function,
satisfying $v(\emptyset)=0$.

If $f_n:\ \Omega\to\reals$, $n\in\nintegers$, and 
$f:\ \Omega\to\reals$ are measurable and Choquet integrable
functions such that $f_n$ are pointwise non-decreasing in $n$ and
$\dsp f=\limf{n}f_n$, $v$-almost surely,
then $\dsp\limf{n}v(f_n)=v(f)$ holds.
\DDD
\thme
\prfb
First we assume that $f_1$ is non-negative valued. Then
$f_n$\,, $n\in\nintegers$, and $f$ are also non-negative.
For $z\in\reals$, put
$\dsp C_n(z)=\{\omega\in\Omega\mid f_n(\omega)>z\}$, $n\in\nintegers$,
and $\dsp {\cal C}(z)=\{C_n(z)\mid n\in\nintegers\}$.
By assumption that $f_n$ is pointwise non-dereasing in $n$,
the limit $\dsp\limf{n}f_n:\ \Omega\to\reals\cup\{+\infty\}$ exists
if we conventionally allow $+\infty$ for a limit value. 
Furthermore it holds that
$\dsp C_1(z)\subset C_2(z)\subset\cdots$ 
and $\dsp \bigcup_{n\in\nintegers} C_n(z)
=\{\omega\in\Omega\mid\limf{n}f_n(\omega)>z\}$.
Define the non-negative functions $g_n:\ \reals\to\preals$, $n\in\nintegers$,
and $g:\ \reals\to\preals$ by
$\dsp g_n(z)=v(C_n(z))$, $n\in\nintegers$, and 
$\dsp g(z)=v(\bigcup_{n\in\nintegers} C_n(z))$.
Since by assumption $v$ is non-decreasing $g_n$ is pointwise non-decreasing
in $n$ and $g_n(z)\le g(z)$, $z\in\reals$.

By assumption $\sgcn\ne\emptyset$ there exists ${\cal I}\in\sgcn$,
hence \lemu{ctblechaininsertion} implies $\dsp {\cal I}_{{\cal C}(z)}\in\sgcn$.
Therefore, assumption of continuity of $v$ implies that
there exists a measure $\dsp\mu_{v,{\cal I}_{{\cal C}}(z)}:\ {\cal F}\to\preals$
defined by \eqnu{extremept} with $\dsp {\cal I}={\cal I}_{{\cal C}(z)}$\,.
It therefore follows from continuity of measure that
\[\arrb{l}\dsp
\limf{n} g_n(z)=\limf{n} v(C_n(z))
=\limf{n} \mu_{v,{\cal I}_{{\cal C}}(z)}(C_n(z))
= \mu_{v,{\cal I}_{{\cal C}}(z)}(\bigcup_{n\in\nintegers} C_n(z))
= v(\bigcup_{n\in\nintegers} C_n(z))
\\ \dsp \phantom{\limf{n} g_n(z)}
=g(z),\ z\in\reals.
\arre\]
For $z\in\reals$, put
$\dsp C(z)=\{\omega\in\Omega\mid f(\omega)>z\}$ and
$\dsp N=\{\omega\in\Omega\mid f(\omega)\ne \limf{n}f_n(\omega)\}$,
Then $\dsp \bigcup_{n\in\nintegers}C_n(z)\cap N^c
\subset C(z)\subset \bigcup_{n\in\nintegers}C_n(z)\cup N$ hols.
Assumption of $v$-almost sure onvergence implies that $N$ is a $v$-null set,
hence we have $\dsp v(C(z))=v(\bigcup_{n\in\nintegers}C_n(z))=g(z)$.

Choquet integrability assumptions of $f_n$ and $f$, 
with \eqnu{Choquetbd} with $x=0$ then imply
\[
v(f)
=\int_{0}^{\infty} g(z)\,dz,
\ \ \ 
v(f_n)
=\int_{0}^{\infty} g_n(z)\,dz,\ n\in\nintegers.
\]
Hence the standard measure theoretic
monotone convergence theorem for $1$-dimensional Lebesgue measure
implies $\dsp \limf{n} v(f_n)=\limf{n} \int_{0}^{\infty} g_n(z)\,dz
=\int_{0}^{\infty} g(z)\,dz=v(f)$,
which proves the claim when the functions are non-negative.

Next, we consider the case that the functions are bounded from below, 
and assume that there exists $\dsp y_0\in\reals$ such that
$\dsp f_1(\omega)\ge y_0$\, $\omega\in\Omega$. By pointwise non-decreasing 
properties of the functions, it follows that
$\dsp f_n(\omega)\ge y_0$\, $n\in\nintegers$, and $f(\omega)\ge y_0$,
for all $\omega\in\Omega$.
Subtracting the constant $y_0$ from the functions, put
$\dsp \tilde{f}_n=f_n-y_0$\,, $n\in\nintegers$, and
$\dsp \tilde{f}=f-y_0$\,.
Then $\dsp \tilde{f}_n$\,, $n\in\nintegers$, and $\tilde{f}$ are
non-negative vauled functions satisfying the assumptions of the Theorem,
hence the proof in the previous paragraphs implies
$\dsp\limf{n}v(\tilde{f}_n)=v(\tilde{f})$, $v$-almost surely.
Also \eqnu{Choquetbd} implies, 
with change of integration variable $z'=z-y_0$\,,
\[\arrb{l}\dsp
v(f)=y_0\,v(\Omega)
+\int_{y_0}^{\infty} v(\{\omega\in\Omega\mid f(\omega)>z\})\,dz
\\ \dsp \phantom{v(f)}
=y_0\,v(\Omega)
+\int_{0}^{\infty} v(\{\omega\in\Omega\mid \tilde{f}(\omega)>z'\})\,dz'
=y_0\,v(\Omega)+v(\tilde{f}),
\arre\]
and similar formula for $f_n$\,, $n\in\nintegers$.
Therefore we have $\dsp\limf{n}v(f_n)=v(f)$, $v$-almost surely.

Finally we consider the general case.
Since $v$ is non-decreasing, \eqnu{Choquetgen} implies that
for any $\eps>0$ there exists $y_0\in\reals$ such that
\begin{equation}
\eqna{monconvthmprf1}
v(f_1)\le y_0\,v(\Omega)
+\int_{y_0}^{\infty} v(\{\omega\in\Omega\mid f_1(\omega)>z\})\,dz 
\le v(f_1)+\eps.
\end{equation}
Round up the values of the functions smaller than $y_0$ and 
put $\dsp \tilde{f}_n=f_n\vee y_0$\,, $n\in\nintegers$, and
$\dsp \tilde{f}=f\vee y_0$\,, where $a\vee b=a$ if $a\ge b$ and
otherwise $a\vee b=b$.
Then $\dsp \tilde{f}_n$\,, $n\in\nintegers$, and $\tilde{f}$ are
bounded from below and satisfy the assumptions of the Theorem,
hence the proof in the previous paragraph implies
\begin{equation}
\eqna{monconvthmprf2}
\limf{n}v(\tilde{f}_n)=v(\tilde{f}),
\end{equation}
$v$-almost surely.
Also $\dsp \tilde{f}_n=f_n\vee y_0$ and \eqnu{Choquetbd} imply
\begin{equation}
\eqna{monconvthmprf3}
\arrb{l}\dsp
v(\tilde{f}_n)=y_0\,v(\Omega)
+\int_{y_0}^{\infty} v(\{\omega\in\Omega\mid \tilde{f}_n(\omega)>z\})\,dz
\\ \dsp \phantom{v(\tilde{f}_n)}
=y_0\,v(\Omega)
+\int_{y_0}^{\infty} v(\{\omega\in\Omega\mid f_n(\omega)>z\})\,dz,
\ n\in\nintegers,
\arre
\end{equation}
hence 
\eqnu{monconvthmprf1} implies
\begin{equation}
\eqna{monconvthmprf4}
0\le  v(\tilde{f}_1)-v(f_1) \le\eps.
\end{equation}
Also, \eqnu{monconvthmprf3} and \eqnu{Choquetgen} imply
\begin{equation}
\eqna{monconvthmprf5}
v(\tilde{f}_n)-v(f_n)
=\lim_{y\to-\infty} \biggl((y_0-y)\,v(\Omega)
-\int_{y}^{y_0} v(\{\omega\in\Omega\mid f_n(\omega)>z\})\,dz
\biggr),\ n\in\nintegers,
\end{equation}
Since $f_n$ is pointwise non-decreasing in $n$ and since $v$ is non-decreasing
$v(\{\omega\in\Omega\mid f_n(\omega)>z\})$ is also non-decreasing in $n$ and
is no larger than $v(\Omega)$. Therefore \eqnu{monconvthmprf5} further implies
\begin{equation}
\eqna{monconvthmprf6}
v(\tilde{f}_1)-v(f_1)\ge v(\tilde{f}_2)-v(f_2)\ge \cdots \ge
\limf{n}(v(\tilde{f}_n)-v(f_n))\ge 0
\end{equation}
The argument between \eqnu{monconvthmprf3} and \eqnu{monconvthmprf6},
applied to the function $f$ and $\tilde{f}$ implies
\begin{equation}
\eqna{monconvthmprf7}
v(\tilde{f})-v(f)
=\lim_{y\to-\infty} \biggl((y_0-y)\,v(\Omega)
-\int_{y}^{y_0} v(\{\omega\in\Omega\mid f(\omega)>z\})\,dz
\biggr)\ge0.
\end{equation}
Combining \eqnu{monconvthmprf4}, \eqnu{monconvthmprf6},
\eqnu{monconvthmprf2}, and \eqnu{monconvthmprf7}, we therefore have
\begin{equation}
\eqna{monconvthmprf8}
\eps \ge 
\limf{n}v(\tilde{f}_n) - \limf{n}v(f_n)=v(\tilde{f})- \limf{n}v(f_n)
\ge v(f)- \limf{n}v(f_n).
\end{equation}
Since by assumption $f_n$ is pointwise non-decreasing in $n$ and
converges $v$-almost surely to $f$,
$v(f_n)\le v(f)$, 
hence $\dsp \limf{n}v(f_n)\le v(f)$.
Since $\eps>0$ is arbitrary $\dsp \limf{n}v(f_n)= v(f)$ follows.
\QED
\prfe
As in standard measure theory, a simple function approximation of 
a non-negative Choquet integrable function $f:\ \Omega\to\preals$ defined by
\begin{equation}
\eqna{simplfcnapp}
f_k:=\sum_{i=1}^{k\,2^k}2^{-k}(i-1)\chrfcn{2^{-k}(i-1)<f\le 2^{-k}i}
+k\chrfcn{\{\omega\in\Omega\mid f(\omega)>k\}},
\ \ k\in\nintegers,
\end{equation}
is an example of non-decreasing sequence of functions which
converges pointwise to $f$ as $k\to\infty$.
\thmu{monconvthm} then implies $\dsp v(f)=\limf{k} v(f_k)$.

\subsection{Representation formula for Choquet integration and coherent risk measure}

The representation \eqnu{submodextremeptrep} for submodular functions
in the framework of \thmu{totorder2modularvarrepbigcoreequiv} implies
corresponding formula for Choquet integrations.
\thmb
\thma{genDelbaenKusuoka}
Let $(\Omega,{\cal F})$ be a measurable space satisfying $\sgcn\ne\emptyset$.
Let $v:\ {\cal F}\to\reals$ be a non-decreasing and continuous submodular
function satisfying $v(\emptyset)=0$.
Then for a bounded measurable function $\dsp f:\ \Omega\to\reals$
\begin{equation}
\eqna{genDelbaenKusuoka}
v(f)=\sup_{{\cal I}\in \sgcn}
 \int_{\Omega} f(\omega)\,\mu_{v,{\cal I}}(d\omega)
\end{equation}
holds.
\DDD
\thme
\prfb
Since $f$ is bounded and measurable,
$f$ is both $v$-integrable and integrable with respect to 
(finite) measures.
Also 
since \eqnu{Choquetbd} implies $\dsp v(f+a)=av(\Omega)+v(f)$,
we may assume from beginning that $f$ is non-negative in the proof of
\eqnu{genDelbaenKusuoka} by adding a positive constant to $f$.
Let $f_k$\,, $k\in\nintegers$, be the series of 
simple function approximations \eqnu{simplfcnapp} of $f$.
For simplicity of notation we follow that of \eqnu{simplefcnChoquetrep}
and write \eqnu{simplfcnapp} as
\[\arrb{l}\dsp
f_k=\sum_{i=1}^{n_k}a^{(k)}_i\chrfcn{A^{(k)}_i}\,,\ \ 
a^{(k)}_i>0,\ i=1,\ldots,n_k-1,\ a^{(k)}_{n_k}=0,\\ \dsp
A^{(k)}_1\subset A^{(k)}_2\subset\cdots\subset A^{(k)}_{n_k}=\Omega,\ \ 
{\cal A}_k=\{A^{(k)}_i\mid i=1,\ldots,n_k\}.
\arre\]
Note that we can put $a^{(k)}_{n_k}=0$ because we assume that $f$ is
non-negative.
Then \eqnu{simplefcnChoquetint} implies 
$\dsp v(f_k)=\sum_{i=1}^{n_k} a^{(k)}_i\,v(A^{(k)}_i)$.

As noted below \eqnu{finiteadditivemeasure},
the definition of the continuity of $v$ implies 
that \eqnu{extremept} uniquely defines a measure
$\dsp \mu_{v,{\cal I}}:\ {\cal F}\to\reals$, for each ${\cal I}\in\sgcn$.
In particular, for $\dsp {\cal I}_{{\cal A}_k}\in\sgcn$\,, the insertion of
a chain $\dsp {\cal A}_k$ to ${\cal I}\in\sgcn$ we have, from
\eqnu{extremept},
\begin{equation}
\eqna{genDelbaenKusuokaprf1}
v(f_k)=\sum_{i=1}^{n_k}a^{(k)}_iv(A^{(k)}_i)
=\sum_{i=1}^{n_k}a^{(k)}_i\mu_{v,{\cal I}_{{\cal A}_k}}(A^{(k)}_i)
=\int_{\Omega} f_k(\omega)\,\mu_{v,{\cal I}_{{\cal A}_k}}(d\omega).
\end{equation}
Also, \lemu{terminalmeasureinC} implies, as in the proof of (a) to (b) in
\thmu{totorder2modularvarrepbigcoreequiv},
$\dsp \mu_{v,{\cal I}}(A)\le v(A)$, $A\in{\cal F}$, ${\cal I}\in\sgcn$.
Therefore,
\begin{equation}
\eqna{genDelbaenKusuokaprf2}
v(f_k)=\sum_{i=1}^{n_k}a^{(k)}_iv(A^{(k)}_i)
\ge\sum_{i=1}^{n_k}a^{(k)}_i\mu_{v,{\cal I}}(A^{(k)}_i)
=\int_{\Omega} f_k(\omega)\,\mu_{v,{\cal I}}(d\omega).
\ \ {\cal I}\in\sgcn.
\end{equation}

Since $f_k$ is pointwise non-decreasing in $k$, 
\begin{equation}
\eqna{genDelbaenKusuokaprf3}
\int_{\Omega} f_k(\omega)\,\mu_{v,{\cal I}}(d\omega)
\le\int_{\Omega} f(\omega)\,\mu_{v,{\cal I}}(d\omega),\ {\cal I}\in\sgcn,
\mbox{ and } v(f_k)\le v(f),\ \ k\in\nintegers.
\end{equation}
The standard monotone convergence theorem for measures implies
\begin{equation}
\eqna{genDelbaenKusuokaprf4}
(\forall \eps>0)(\forall {\cal I}\in\sgcn) \exists k_0\in\nintegers;\ 
(\forall k\ge k_0)\ 
\int_{\Omega} f(\omega)\,\mu_{v,{\cal I}}(d\omega)\le
\int_{\Omega} f_k(\omega)\,\mu_{v,{\cal I}}(d\omega)+\eps
\end{equation}
and \thmu{monconvthm} implies
\begin{equation}
\eqna{genDelbaenKusuokaprf5}
(\forall \eps>0) \exists k_0\in\nintegers;\ (\forall k\ge k_0)\ 
v(f)\le v(f_k)+\eps.
\end{equation}

Combining \eqnu{genDelbaenKusuokaprf2}, \eqnu{genDelbaenKusuokaprf3}, and
\eqnu{genDelbaenKusuokaprf4}, we see that 
$\dsp v(f)\ge 
\int_{\Omega} f(\omega)\,\mu_{v,{\cal I}}(d\omega)-\eps$
for any ${\cal I}\in\sgcn$ and $\eps>0$. Hence
\begin{equation}
\eqna{genDelbaenKusuokaprf6}
v(f)\ge \sup_{{\cal I}\in\sgcn}
\int_{\Omega} f(\omega)\,\mu_{v,{\cal I}}(d\omega).
\end{equation}
Similarly, 
combining \eqnu{genDelbaenKusuokaprf1}, \eqnu{genDelbaenKusuokaprf3}, and
\eqnu{genDelbaenKusuokaprf5}, we see that for any $\eps>0$ there exists
$k_0\in\nintegers$ such that if $k\ge k_0$ then
$\dsp v(f)\le 
\int_{\Omega} f(\omega)\,\mu_{v,{\cal I}_{{\cal A}_k}}(d\omega)+\eps
\le
\sup_{{\cal I}\in\sgcn}\int_{\Omega} f(\omega)\,\mu_{v,{\cal I}}(d\omega)+\eps
$.
Hence $\dsp v(f)\le \sup_{{\cal I}\in\sgcn}
\int_{\Omega} f(\omega)\,\mu_{v,{\cal I}}(d\omega)$, which, with
\eqnu{genDelbaenKusuokaprf6}, implies \eqnu{genDelbaenKusuoka}.
\QED
\prfe
\corb
\cora{Choquet2cohrentriskmeasure}
Let $(\Omega,{\cal F})$ be a measurable space satisfying $\sgcn\ne\emptyset$.
Let $v:\ {\cal F}\to\reals$ be a non-decreasing and continuous submodular
function satisfying $v(\emptyset)=0$, and
define a function $\rho$ on a set of bounded measurable functions 
by $\dsp \rho(f)=\frac{v(-f)}{v(\Omega)}$.
Then the following hold.
\ittb
\item[non-negativity: ]
If $\dsp f\ge 0$ then $\dsp \rho(f)\le 0$,
\item[subadditivity: ]
$\dsp \rho(f+g)\le \rho(f)+\rho(g)$,
\item[positive homogeneity: ]
If $\dsp \lambda\ge0$ then $\dsp\rho(\lambda f)=\lambda \rho(f)$,
\item[translational invariance: ]
If $a\in\reals$ then $\dsp \rho(f+a)=\rho(f)-a$.
\DDD
\itte
\core
\prfb
The claims are straightforward consequences of
\eqnu{genDelbaenKusuoka} in \thmu{genDelbaenKusuoka}.
For example,
\[\arrb{l}\dsp
\rho(f+g)=-\frac1{v(\Omega)}\,\inf_{{\cal I}\in \sgcn}
 \int_{\Omega} (f(\omega)+g(\omega))\,\mu_{v,{\cal I}}(d\omega)
\\ \dsp \phantom{\rho(f+g)}
\le -\frac1{v(\Omega)}\,\inf_{{\cal I}\in \sgcn}
 \int_{\Omega} f(\omega)\,\mu_{v,{\cal I}}(d\omega)
-\frac1{v(\Omega)}\,\inf_{{\cal I}\in \sgcn}
 \int_{\Omega} g(\omega))\,\mu_{v,{\cal I}}(d\omega)
=\rho(f)+\rho(g),
\arre\]
which proves subadditivity.
Proofs of non-negativity, positive homogeneity, and translational invariance
are easier.
\QED
\prfe
In the field of mathematical finance,
the set of properties in \coru{Choquet2cohrentriskmeasure} is known
to be the definition that $\rho$ is a coherent risk measure
\cite{risk99,risk01,riskmSK,19Peng}.
Thus \coru{Choquet2cohrentriskmeasure} relates
the Choquet integration with respect to a sub-modular set function
to a coherent risk measure.
This observation motivates reformulating in our framework the results
found in mathematical finance, which we consider in the next section
\secu{lawinvcrm}.

\section{Choquet integration for uniform case and law invariant coherent risk measure}
\seca{lawinvcrm}

In \cite{riskmSK}, coherent risk measure is studied on probability spaces. 
Given a probability measure $\prb{}$, the law invariance of a
set function $v$ means that 
$\prb{A}=\prb{B}$ implies $v(A)=v(B)$.
In our approach, we started with a measurable space $(\Omega,{\cal F})$ 
without probability measure, hence 
we are free to intruduce a (finite) measure $\nu:\ {\cal F}\to\reals$ and
define law invariance as $\nu(A)=\nu(B)$ implies $v(A)=v(B)$.

In addition, it is assumed in \cite{riskmSK} that $(\Omega,{\cal F}, \prb{})$
is a standard probability space and that $\prb{}$ is non-atomic,
and with these assumptions, the proof is reduced to 
the $1$-dimensional Borel measurable space on an interval
$(\Omega,{\cal F})=([0,1),{\cal B}_1([0,1)))$, and
$\prb{}$ specified as the $1$-dimensional Lebesgue measure.
Here, we keep the only assumption $\sgcn\ne\emptyset$,
existence of a chain generating ${\cal F}$, for $(\Omega,{\cal F})$,
and see how a formula corresponding to that studied in \cite{riskmSK} 
is deduced in our framework.

As a simple example of submodular function whose dependence on the
variable $A\in{\cal F}$ is given through $\nu(A)$, we note the follwing.
In the following, for $a,b\in\reals$ such that $a\le b$ we write 
$a\wedge b=a$ and $a\vee b=b$.
For example, $a+b=(a\vee b)+(a\wedge b)$ and $a\vee b\ge a\wedge b$ hold.
\prpb
\prpa{uniformsubmodularfcn}
Let $(\Omega,{\cal F})$ be a measurable space, $\nu:\ {\cal F}\to\reals$
a measure on the space, and $c\in\reals$.
Then the set function $v:\ {\cal F}\to\reals$ defined by
$\dsp v(A)=c\wedge \nu(A)$, $A\in{\cal F}$, is submodular, i.e., satisfies
\eqnu{submodular}.
\DDD
\prpe
\prfb
Let $A,B\in{\cal F}$. Then additivity, non-negativity, and monotonicity 
of the measure $\nu$ imply
\[ \arrb{l}\dsp v(A)+v(B)-v(A\cup B)-v(A\cap B)
\\ \dsp {}
=\left\{\arrb{ll} \dsp 
c+c-c-c=0, & \dsp c\le \nu(A\cap B), \\ \dsp
c+c-\nu(A\cap B)-c\ge 0, &\dsp \nu(A\cap B)< c\le \nu(A)\wedge \nu(B), \\ \dsp
(\nu(A)\wedge \nu(B))+c-\nu(A\cap B)-c\ge0,
&\dsp  \nu(A)\wedge \nu(B)< c\le  \nu(A)\vee \nu(B), \\ \dsp
\nu(A)+\nu(B)-\nu(A\cap B)-c & \\ \dsp
=\nu(A\cup B)-c\ge 0,
&\dsp  \nu(A)\vee \nu(B)< c\le \nu(A\cup B), \\ \dsp
\nu(A)+\nu(B)-\nu(A\cup B)-\nu(A\cap B)=0, &\dsp \nu(A\cup B)\le c.
\arre\right.
\arre\]
which proves \eqnu{submodular}.
\QED
\prfe
Submodularlity \eqnu{submodular} is preserved by summation and
multiplication of positive reals. This leads to considering $v$ of a form
$\dsp v(A)=\int_{\preals} g(z)\wedge \nu(A)\,dz$ for some 
non-negative function $g$.
If we choose the parameter $z$ so that $g$ is monotone,
we could consider $g$ as a distribution function of a measure.
These considerations suggest considering a following form for a
submodular set function $v$ whose depedence on variable $A\in{\cal F}$
enters through $\nu(A)$.
\thmb
\thma{uniformsubmodularfcn}
Let $(\Omega,{\cal F})$ be a measurable space satisfying $\sgcn\ne\emptyset$,
and $\nu\in {\cal M}(\Omega)$.
For a finite measure $\mu\in {\cal M}(\Omega)$ satisfying $\mu\ll\nu$,
denote the Radon--Nykodim derivative by 
$\dsp \diff{\mu}{\nu}:\ \Omega\to\preals$\,,
and define the distribution function $F_{\mu}:\ \preals\to\preals$ by
\begin{equation}
\eqna{distributionfcn}
F_{\mu}(y)
=\nu(\{\omega\in\Omega\mid \diff{\mu}{\nu}(\omega)\le y\}),\ y\in\preals,
\end{equation}
and define a set function $v_{\mu}:\ {\cal F}\to\reals$ by
\begin{equation}
\eqna{fcnintgl2submodular}
v_{\mu}(A)
=\int_0^{\infty} (\nu(\Omega)-F_{\mu}(z))\wedge \nu(A)\,dz,\ A\in{\cal F}.
\end{equation}
Then the following hold.
\itmb
\item
$v_{\mu}$ is non-decreasing, continuous, submodular, and
satisfies $v_{\mu}(\emptyset)=0$.

\item
It holds that
\begin{equation}
\eqna{fmuincreasingsaddgencompliment}
v_{\mu}(\{\omega\in\Omega\mid \diff{\mu}{\nu}(\omega)>y\})
=\mu(\{\omega\in\Omega\mid \diff{\mu}{\nu}(\omega)>y\}),\ y\in\preals.
\end{equation}

\item
For a non-negative valued measurable function $f:\ \Omega\to\preals$,
let $\dsp v_{\mu}(f)$ denote the Choquet integration of $f$ with respect to
$\dsp v_{\mu}$ in \eqnu{Choquetbd}. Then
\begin{equation}
\eqna{genDKChoquetsubmodular}
v_{\mu}(f)
\ge\sup\{ \int_{\Omega} f(\omega) \mu'(d\omega)\mid 
\mu'\in{\cal M}(\Omega),\ \mu'\ll \nu,\ F_{\mu'}=F_{\mu}\},
\end{equation}
holds.
In particular, it holds that
\begin{equation}
\eqna{DelbaenKusuokabd}
v_{\mu}(A)\ge
\sup\{ \mu'(A)\mid \mu'\in{\cal M}(\Omega),\ \mu'\ll \nu,\ F_{\mu'}=F_{\mu}\},
\ A\in{\cal F}.
\end{equation}
\DDD
\itme
\thme
\prfb
\itmb
\item
Since $\nu$ is a measure, $\nu(\emptyset)=0$, 
hence $\dsp v_{\mu}(\emptyset)=0$.
Non-decreasing property is obvious from the same property for measures and 
that integration preserves inequality.
To prove that $\dsp v_{\mu}$ is submodular,
substitute $c$ in \prpu{uniformsubmodularfcn} 
with $\dsp \nu(\Omega)-F_{\mu}(z)$ and integrate over $z$,
we see from \prpu{uniformsubmodularfcn} that $\dsp v_{\mu}$ satisfies
\eqnu{submodular}.

To prove continuity, since we have shown that $\dsp v_{\mu}$ is
non-decreasing and submodular,
\prpu{conti24equivconti23onsubmod} now implies that it suffices to prove
\eqnu{upperconti} and \eqnu{lowerconti} for $\dsp v=v_{\mu}$\,.
Since $\nu$ is a measure, 
\eqnu{upperconti} and \eqnu{lowerconti} hold for $\dsp v=\nu$.
This and monotone convergence theorem applied to \eqnu{fcnintgl2submodular}
imply \eqnu{upperconti} and \eqnu{lowerconti} for $\dsp v=v_{\mu}$\,.

\item
Note a basic formula
\begin{equation}
\eqna{Fubininonsigmafinite}
\int_0^{\infty}\nu(f>z)\,dz
:=\int_0^{\infty}\nu(\{\omega\in\Omega\mid f(\omega)>z\})\,dz
=\int_{\Omega}f\,d\nu,
\end{equation}
valid for any 
non-negative measurable function $f:\ \Omega\to\preals$.
In particular, for $v_{\mu}$ as in \eqnu{fcnintgl2submodular}, we have
$\dsp v_{\mu}(\Omega)=\int_{\Omega} \diff{\mu}{\nu}d\nu=\mu(\Omega)$.
By replacement 
$f\mapsto f\vee a$,
\eqnu{Fubininonsigmafinite} also implies
that for any non-negative measurable function 
$f:\ \Omega\to\preals$ and a non-negative constant $a\ge0$,
\begin{equation}
\eqna{Fubininonsigmafiniteindefinite2}
\arrb{l}\dsp
\int_a^{\infty} \nu(f>z)\,dz
=\int_0^{\infty} \nu(f\vee a>z)\,dz-a\nu(\Omega)
=\int_{\Omega}f\vee a\,d\nu-a\nu(\Omega)
\\ \dsp\phantom{\int_a^{\infty} \nu(f>z)\,dz}
=\int_{f\ge a}f\,d\nu-a\nu(f\ge a) =\int_{f> a}f\,d\nu-a\nu(f> a).
\arre
\end{equation}

Using \eqnu{Fubininonsigmafiniteindefinite2} with
$\dsp f=\diff{\mu}{\nu}$ and $\dsp a=y$ in \eqnu{fcnintgl2submodular} we have
\[\arrb{l}\dsp
v_{\mu}(\{\omega\in\Omega\mid \diff{\mu}{\nu}(\omega)>y\})
\\ \dsp {}
=\int_y^{\infty}\nu((\{\omega\in\Omega\mid \diff{\mu}{\nu}(\omega)>z\})\,dz
+y\,\nu((\{\omega\in\Omega\mid \diff{\mu}{\nu}(\omega)>y\})
\\ \dsp \phantom{}
=\int_{\dsp\diff{\mu}{\nu}(\omega)>y}\diff{\mu}{\nu}(\omega)\,\nu(d\omega)
=\mu(\{\omega\in\Omega\mid \diff{\mu}{\nu}(\omega)>y\}
\arre\]
which proves \eqnu{fmuincreasingsaddgencompliment}.

\item
Let $\mu':\ {\cal F}\to\preals$ be a measure satisfying 
$\mu'\ll \nu$ and $\dsp F_{\mu'}=F_{\mu}$\,.
The definition \eqnu{Choquetbd} of Choquet integration,
an elementary inequality
$\dsp \nu(A)\wedge\nu(B)\ge\nu(A\cap B)$, valid for any $A,B\in{\cal F}$,
and \eqnu{Fubininonsigmafinite} then imply
\[\arrb{l}\dsp
v_{\mu}(f)=\int_{0}^{\infty} v_{\mu}(\{\omega\in\Omega\mid f(\omega)>z\})\,dz
\\ \dsp \phantom{v_{\mu}(f)}
=\int_{0}^{\infty} 
\int_0^{\infty} (\nu(\Omega)-F_{\mu}(y))\wedge
\nu(\{\omega\in\Omega\mid f(\omega)>z\})\,dy\,dz
\\ \dsp \phantom{v_{\mu}(f)}
=\int_{0}^{\infty} 
\int_0^{\infty} (\nu(\Omega)-F_{\mu'}(y))\wedge
\nu(\{\omega\in\Omega\mid f(\omega)>z\})\,dy\,dz
\\ \dsp \phantom{v_{\mu}(f)}
=\int_{0}^{\infty} 
\int_0^{\infty} (\nu(\{\omega\in\Omega\mid \diff{\mu'}{\nu}(\omega)>y\}))\wedge
\nu(\{\omega\in\Omega\mid f(\omega)>z\})\,dy\,dz
\\ \dsp \phantom{v_{\mu'}(f)}
\ge\int_{0}^{\infty} 
\int_0^{\infty} (\nu(\{\omega\in\Omega\mid \diff{\mu'}{\nu}(\omega)>y\})
\cap\{\omega\in\Omega\mid f(\omega)>z\})\,dy\,dz
\\ \dsp \phantom{v_{\mu'}(f)}
=\int_{0}^{\infty} \int_{\Omega}
\int_0^{\infty} \chrfcn{\diff{\mu'}{\nu}(\omega)>y}
\chrfcn{f(\omega)>z}\,dy\,\nu(d\omega)\,dz
\\ \dsp \phantom{v_{\mu'}(f)}
=\int_{0}^{\infty} \int_{\Omega}
\chrfcn{f(\omega)>z}\,\diff{\mu'}{\nu}(\omega)\,\nu(d\omega)\,dz
\\ \dsp \phantom{v_{\mu'}(f)}
=\int_{0}^{\infty} \int_{\Omega}
\chrfcn{f(\omega)>z}\,\mu'(d\omega)\,dz
=\int_{0}^{\infty}\mu'(\{\omega\in\Omega\mid f(\omega)>z\})\,dz
\\ \dsp \phantom{v_{\mu'}(f)}
=\int_{\Omega}f(\omega)\,\mu'(d\omega),
\arre\]
which proves \eqnu{genDKChoquetsubmodular}.

By choosing $\dsp f=\chrfcn{A}$ in \eqnu{genDKChoquetsubmodular},
\eqnu{DelbaenKusuokabd} follows.
\QED
\itme
\prfe

As in \cite[Theorem~7]{riskmSK}, we assume comonotonicity for
further results.
We say that the functions $f:\ \Omega\to\reals$ and $g:\ \Omega\to\reals$
are comonotone if
$\dsp 
 \{ \{\omega\in\Omega\mid f(\omega)\le z\}\mid z\in\reals \}\cup
 \{ \{\omega\in\Omega\mid g(\omega)\le z\}\mid z\in\reals \}
$
is a chain, i.e.,
for each $\dsp (y,z)\in\nreals{2}$, either
$\dsp \{\omega\in\Omega\mid g(\omega)\le y\}
\supset \{\omega\in\Omega\mid f(\omega)\le z\}$
or
$\dsp \{\omega\in\Omega\mid g(\omega)\le y\}
\subset \{\omega\in\Omega\mid f(\omega)\le z\}$ holds.
\thmb
\thma{comonotoneChoquetExprep}
Let $(\Omega,{\cal F})$ be a measurable space satisfying $\sgcn\ne\emptyset$,
and $\nu\in {\cal M}(\Omega)$.
For a finite measure $\mu\in {\cal M}(\Omega)$ satisfying $\mu\ll\nu$,
define a set function $v_{\mu}:\ {\cal F}\to\reals$ by 
\eqnu{fcnintgl2submodular}.
For $y\ge0$ 
put $\dsp I_{\mu,y}=\{\omega\in\Omega\mid \diff{\mu}{\nu}(\omega)> y\}$ and
put $\dsp {\cal I}_{\mu}=\{\emptyset,\Omega\}\cup 
\{I_{\mu,y}\mid  y\in\preals\}$ and assume that
\begin{equation}
\eqna{mugeneratedsadd}
{\cal I}_{\mu}\in\sgcn \mbox{ and }
\img F_{\mu}\supset\img\nu,
\end{equation}
where $\img$ denotes the image of the map.

Then if a non-negative valued measurable function $f:\ \Omega\to\preals$
is comonotone with $\dsp \diff{\mu}{\nu}$,
the equality in \eqnu{genDKChoquetsubmodular} is attained by $\mu$.
Namely,
\begin{equation}
\eqna{DelbaenKusuokaf2}
v_{\mu}(f)=\int_{\Omega} f(\omega)\, \mu(d\omega)
=\sup\{ \int_{\Omega} f(\omega) \mu'(d\omega)\mid 
\mu'\in{\cal M}(\Omega),\ \mu'\ll \nu,\ F_{\mu'}=F_{\mu}\}
\end{equation}
holds.
\DDD
\thme
\prfb
Note first that
by assumptions there exists non-decreasing function
$\dsp y_0:\ \preals\to\preals$ such that
\begin{equation}
\eqna{comonotoneChoquetExprepprf1}
\arrb{l}\dsp
\nu(\{\omega\in\Omega\mid f(\omega)>y\})
=\nu(\Omega)-F_{\mu}(y_0(y))
\\ \dsp {}
=\nu(\{\omega\in\Omega\mid \diff{\mu}{\nu}(\omega)> y_0(y)\}),
\ \ y\in\preals.
\arre
\end{equation}
The definiton \eqnu{fcnintgl2submodular} of $v_{\mu}$\,,
\eqnu{comonotoneChoquetExprepprf1}, and
\eqnu{fmuincreasingsaddgencompliment} in \thmu{uniformsubmodularfcn}
imply
\begin{equation}
\eqna{comonotoneChoquetExprepprf2}
\arrb{l}\dsp
v_{\mu}(\{\omega\in\Omega\mid f(\omega)>y\})
=\int_0^{\infty} (\nu(\Omega)-F_{\mu}(z))\wedge 
\nu(\{\omega\in\Omega\mid f(\omega)>y\})\,dz
\\ \dsp {}
=\int_0^{\infty} (\nu(\Omega)-F_{\mu}(z))\wedge 
\nu(\{\omega\in\Omega\mid \diff{\mu}{\nu}(\omega)> y_0(y)\})\,dz
\\ \dsp {}
= v_{\mu}(\{\omega\in\Omega\mid \diff{\mu}{\nu}(\omega)>y_0(y)\})
\\ \dsp {}
= \mu(\{\omega\in\Omega\mid \diff{\mu}{\nu}(\omega)>y_0(y)\}),
\ y\in\preals.
\arre
\end{equation}

Denote the set difference by $\triangle$, 
so that $A\triangle B=(A\cap B^c)\cup (A^c\cap B)$.
By the comonotonicity assumption, either
$\dsp \{\omega\in\Omega\mid f(\omega)>y\}\subset 
\{\omega\in\Omega\mid \diff{\mu}{\nu}(\omega)> y_0(y)\}$
or
$\dsp \{\omega\in\Omega\mid f(\omega)>y\}\supset 
\{\omega\in\Omega\mid \diff{\mu}{\nu}(\omega)> y_0(y)\}$,
and in either case \eqnu{comonotoneChoquetExprepprf1} implies that
the difference is measure $0$. Hence
\[
\nu(\{\omega\in\Omega\mid f(\omega)>y\}\triangle
\{\omega\in\Omega\mid \diff{\mu}{\nu}(\omega)>y_0(y)\})=0,
\ \ y\in\preals.
\]
By assumption $\mu\ll\nu$ we then have
$\dsp
\mu(\{\omega\in\Omega\mid f(\omega)>y\}\triangle
\{\omega\in\Omega\mid \diff{\mu}{\mu}(\omega)>y_0(y)\})=0$,
which further implies
$\dsp
\mu(\{\omega\in\Omega\mid f(\omega)>y\})
=\mu(\{\omega\in\Omega\mid \diff{\mu}{\mu}(\omega)>y_0(y)\})$,
$y\in\preals$.
Substituting this in \eqnu{comonotoneChoquetExprepprf2} 
we have
\[
v_{\mu}(\{\omega\in\Omega\mid f(\omega)>y\})
=\mu(\{\omega\in\Omega\mid f(\omega)>y\}),	\ y\in\preals.
\]
Using this in the definition \eqnu{Choquetbd} of Choquet integration
and using \eqnu{Fubininonsigmafinite} we arrive at
\[\arrb{l}\dsp
v_{\mu}(f)=\int_{0}^{\infty} v_{\mu}(\{\omega\in\Omega\mid f(\omega)>y\})\,dy
\\ \dsp \phantom{v_{\mu}(f)}
=\int_{0}^{\infty} \mu(\{\omega\in\Omega\mid f(\omega)>y\})\,dy
=\int_{\Omega} f(\omega)\, \mu(d\omega).
\arre\]
This with \eqnu{genDKChoquetsubmodular} in
\thmu{uniformsubmodularfcn} implies \eqnu{DelbaenKusuokaf2}.
\QED
\prfe

\thmu{comonotoneChoquetExprep} corresponds to \cite[Theorem~7]{riskmSK}
but the notation apparently is quite different.
Before closing this section,
we briefly look into the correspondence of notation 
with the reference, for convenience.
As noted in \coru{Choquet2cohrentriskmeasure} in the previous section,
$v(f)$ corresponds to $\rho(-X)$ in \cite{riskmSK}.
The right continuous inverse of a distribution function $F$
(in the sense of \eqnu{densitydistrfcninv} below) 
is denoted by $Z(x,F)$ in the same reference.
$\dsp v_{\mu}$ has an expression using the inverse of the distribution
function \eqnu{distributionfcn}.
\prpb
\prpa{uniformsubmodularfcn3}
Let $(\Omega,{\cal F})$ be a measurable space satisfying $\sgcn\ne\emptyset$,
and $\nu\in {\cal M}(\Omega)$.
As in \thmu{uniformsubmodularfcn},
let $\mu\in {\cal M}(\Omega)$ satisfying $\mu\ll\nu$ be a finite measure
whose Radon--Nykodim derivative is $\dsp \diff{\mu}{\nu}:\ \Omega\to\preals$\,,
and $F_{\mu}:\ \preals\to\preals$ be its distribution function
\eqnu{distributionfcn},
and $v_{\mu}:\ {\cal F}\to\reals$ defined by \eqnu{fcnintgl2submodular}.
Define $F_{\mu}^{-1}$, the right continuous inverse function of $F_{\mu}$, by
\begin{equation}
\eqna{densitydistrfcninv}
F_{\mu}^{-1}(\beta)= \inf\{ z\in\preals\mid F_{\mu}(z)>\beta\},
\ \beta\in\preals.
\end{equation}
Then $v_{\mu}$ has an expression
\begin{equation}
\eqna{invfcnintgl2submodular}
v_{\mu}(A)=\int_{\nu(A^c)}^{\nu(\Omega)} F_{\mu}^{-1}(\beta)\,d\beta,
\ A\in{\cal F}.
\end{equation}
\DDD
\prpe
\lemb
\lema{uniformsubmodularfcn}
Let $F:\ \preals\to\preals$ be a non-decreasing right continuous 
non-negative valued function on non-negative reals, 
and $F^{-1}:\ \preals\to\preals$
its right continuous inverse function defined by
\begin{equation}
\eqna{fdensitydistrfcninv}
F^{-1}(\beta)= \inf\{ z\in\preals\mid F(z)>\beta\},\ \beta\in\preals\,.
\end{equation}

If $a\in\preals$ and $\alpha\in\preals$ satisfy
$F(a)=\alpha$ or $\dsp F^{-1}(\alpha)=a$ then
\begin{equation}
\eqna{invfcnintglint}
\int_0^{\alpha} F^{-1}(\beta)\,d\beta+\int_0^{a}F(z)\,dz=a\,\alpha
\end{equation}
holds.
If, in addition to $F(a)=\alpha$ or $\dsp F^{-1}(\alpha)=a$,
$\dsp F(+\infty):=\limf{z} F(z)<+\infty$ and
$\dsp \int_0^{\infty} (F(+\infty)-F(z))\,dz<+\infty$ hold,
then
\begin{equation}
\eqna{invfcnintglext}
\int_{\alpha}^{F(+\infty)} F^{-1}(\beta)\,d\beta
-\int_a^{+\infty}(F(+\infty)-F(z))\,dz=a\,(F(+\infty)-\alpha)
\end{equation}
also holds.
\DDD
\leme
\prfb
The formulas hold because both hand sides of each formula
are $2$-dimensional Lebesgue measure (area) of the same rectangle;
both hand sides of \eqnu{invfcnintglint} are area of 
the rectangle $[0,a]\times[0,\alpha]$, and
both hand sides of \eqnu{invfcnintglint} are area of 
the rectangle $[0,a]\times[\alpha,F(+\infty)]$.
\QED
\prfe
%
%
\prfofb{\protect\prpu{uniformsubmodularfcn3}}
\item
Using \eqnu{Fubininonsigmafiniteindefinite2} with
$\dsp f=\diff{\mu}{\nu}$ and $\dsp a=y$ again, we have
\begin{equation}
\eqna{Fubinidensity}
\arrb{l}\dsp
\int_y^{\infty} (\nu(\Omega)-F_{\mu}(z))\,dz
=\int_y^{\infty} \nu(\{\omega\in\Omega\mid \diff{\mu}{\nu}(\omega)>z\})\,dz
\\ \dsp {}
=\mu(\{\omega\in\Omega\mid \diff{\mu}{\nu}(\omega)>y\})
-y\,(\nu(\Omega)-F_{\mu}(y)).
\arre
\end{equation}
In particular, 
$\dsp\int_0^{\infty} (\nu(\Omega)-F_{\mu}(z))\,dz\le\mu(\Omega)<\infty$.

Let $\dsp A\in{\cal F}$, and put
\[\arrb{l}\dsp y_A=\inf\{z\in\preals\mid
 \nu(\{\omega\in\Omega\mid\diff{\mu}{\nu}(\omega)>z\})<\nu(A)\}
\\ \dsp\phantom{y_A}
= \inf\{z\in\preals\mid F_{\mu}(z)>\nu(A^c)\}
=F_{\mu}^{-1}(\nu(A^c)).
\arre\]
Then from \eqnu{fcnintgl2submodular} we have
\[\arrb{l}\dsp
v_{\mu}(A)=\int_0^{\infty} (\nu(\Omega)-F_{\mu}(z))\wedge \nu(A)\,dz
=y_A\,\nu(A)+\int_{y_A}^{\infty} (\nu(\Omega)-F_{\mu}(z))\,dz.
\\ \dsp\phantom{v_{\mu}(A)}
=\nu(A)\,F_{\mu}^{-1}(\nu(A^c))
+\int_{F_{\mu}^{-1}(\nu(A^c))}^{\infty} (\nu(\Omega)-F_{\mu}(z))\,dz.
\arre\]
Substituting $\dsp \alpha =\nu(A^c)$ and $\dsp a=F_{\mu}^{-1}(\nu(A^c))$
in \eqnu{invfcnintglext}, we have
\eqnu{invfcnintgl2submodular}.
\QED
\prfofe
\prpb
\prpa{FubiniChoquetintegration}
Let $(\Omega,{\cal F})$ be a measurable space satisfying $\sgcn\ne\emptyset$,
and $\nu\in {\cal M}(\Omega)$.
For a finite measure $\mu\in {\cal M}(\Omega)$ satisfying $\mu\ll\nu$,
denote the distribution function of the 
Radon--Nykodim derivative $\dsp \diff{\mu}{\nu}$ by 
$F_{\mu}:\ \preals\to\preals$, as in \eqnu{distributionfcn},
and define a set function $v_{\mu}:\ {\cal F}\to\reals$ by
\eqnu{fcnintgl2submodular}.

Then, for a non-negative valued measurable function $f:\ \Omega\to\preals$,
the Choquet integration $v(f)$ of $f$ with respect to $v$ satisfies
\begin{equation}
\eqna{FubiniChoquetintegrationsubmodular}
v_{\mu}(f)=\int_0^{\nu(\Omega)} F_{\mu}^{-1}(\beta)\,F_f^{-1}(\beta)\,d\beta,
\end{equation}
where $F_f$ is the distribution function of $f$ defined by
\begin{equation}
\eqna{distributionfcnf}
F_f(z)=\nu(\{\omega\in\Omega\mid f(\omega)\le z\}),\ z\in\preals,
\end{equation}
and $\dsp F_{\mu}^{-1}$ and $\dsp F_f^{-1}$ are respectively
the right continuous inverse function of $F_{\mu}$ and $F_f$
defined as in \eqnu{densitydistrfcninv}.
\DDD
\prpe
\prfb
The definitions \eqnu{Choquetbd}, \eqnu{fcnintgl2submodular}, and
$\nu(\{\omega\in\Omega\mid f(\omega)>y\})=\nu(\Omega)-F_f(y)$ imply
\begin{equation}
\eqna{FubiniChoquetintegrationprf1}
\arrb{l}\dsp
v_{\mu}(f)
=\int_{0}^{\infty} v_{\mu}(\{\omega\in\Omega\mid f(\omega)>y\})\,dy
\\ \dsp\phantom{v_{\mu}(f)}
=\int_{\preals^2}
(\nu(\Omega)-F_{\mu}(z))\wedge \nu(\{\omega\in\Omega\mid f(\omega)>y\})\,dy\,dz
\\ \dsp\phantom{v_{\mu}(f)}
=\int_{\preals^2}
(\nu(\Omega)-F_{\mu}(z))\wedge (\nu(\Omega)-F_f(y))\,dy\,dz
\\ \dsp\phantom{v_{\mu}(f)}
=\int_{F_f(y)>F_{\mu}(z)} (\nu(\Omega)-F_f(y))\,dy\,dz
+\int_{F_{\mu}(z)\ge F_f(y)} (\nu(\Omega)-F_{\mu}(z))\,dy\,dz.
\arre
\end{equation}
For the first term in the right hand side, 
we perform the $y$ integration first.
Note that since we chose $F_f^{-1}$ to be right continuous,
\[
y>F_f^{-1}(F_{\mu}(z))\ \Rightarrow\ F_f(y)>F_{\mu}(z)\ \Rightarrow\ 
y\ge F_f^{-1}(F_{\mu}(z))\ (\ \Rightarrow\ F_f(y)\ge F_{\mu}(z)\ )
\]
holds.
Changing the integration variable from $z$ to $y$ in \eqnu{invfcnintglext},
choosing $\dsp F=F_f$,  $\dsp \alpha=F_{\mu}(z)$, 
and $\dsp a=F_f^{-1}(F_{\mu}(z))$,
and noting $F_{\mu}(+\infty)=F_f(+\infty)=\nu(\Omega)$ and
recalling that a set of countable points has zero 
$1$-dimensional Lebesgue measure, we have
\[\arrb{l}\dsp
\int_{F_f(y)>F_{\mu}(z)} (\nu(\Omega)-F_f(y))\,dy
=\int_{F_f^{-1}(F_{\mu}(z))}^{+\infty}(\nu(\Omega)-F_f(y))\,dy
\\ \dsp {}
=\int_{F_{\mu}(z)}^{\nu(\Omega)} F_f^{-1}(\beta)\,d\beta
-F_f^{-1}(F_{\mu}(z))\,(\nu(\Omega)-F_{\mu}(z)).
\arre\]
Using this in \eqnu{FubiniChoquetintegrationprf1} with Fubini's theorem and
noting
\[
F_{\mu}(z)<\beta\ \Rightarrow\ z<F^{-1}(\beta)\ \Rightarrow\ F_{\mu}(z)\le\beta
\ (\ \Rightarrow\ z\le F^{-1}(\beta)\ ),
\]
we have
\[\arrb{l}\dsp
v_{\mu}(f)
=\int_{F_{\mu}(z)\le \beta\le \nu(\Omega)} F_f^{-1}(\beta)\,d\beta\,dz
\\ \dsp\phantom{v_{\mu}(f)=}
-\int_0^{\infty} F_f^{-1}(F_{\mu}(z))\,(\nu(\Omega)-F_{\mu}(z))\,dz
+\int_{F_{\mu}(z)\ge F_f(y)} (\nu(\Omega)-F_{\mu}(z))\,dy\,dz
\\ \dsp\phantom{v_{\mu}(f)}
=\int_0^{\nu(\Omega)} F_{\mu}^{-1}(\beta)\,F_f^{-1}(\beta)\,d\beta
\\ \dsp\phantom{v_{\mu}(f)=}
-\int_0^{\infty} F_f^{-1}(F_{\mu}(z))\,(\nu(\Omega)-F_{\mu}(z))\,dz
+\int_{F_f(y)\le F_{\mu}(z)} (\nu(\Omega)-F_{\mu}(z))\,dy\,dz.
\arre\]
Noting that
\[
(\ F_f(y)< F_{\mu}(z)\ \Rightarrow\ )\ y< F_f^{-1}(F_{\mu}(z))
\ \Rightarrow\ F_f(y)\le F_{\mu}(z) \Rightarrow\ y\le F_f^{-1}(F_{\mu}(z)),
\]
we see that the last $2$ terms cancel, which proves
\eqnu{FubiniChoquetintegrationsubmodular}.
\QED
\prfe
%
The expression in the right hand side of 
\eqnu{FubiniChoquetintegrationsubmodular} corresponds to 
the left hand side of the formula in the statement of \cite[Lemma~11]{riskmSK}.

\cite[Theorem~7]{riskmSK} which corresponds to
\thmu{comonotoneChoquetExprep} uses yet another expression.
Let $m$ be the measure on a unit interval $(0,1]$
defined as a Stieltjes measure satisfying
\begin{equation}
\eqna{rskmSKm}
\arrb{l}\dsp
m((0,\gamma])=\frac1{\mu(\Omega)}\int_{
\nu(\{\omega\in\Omega\mid \diff{\mu}{\nu}(\omega)>y\})\le\gamma\nu(\Omega))}
\nu(\{\omega\in\Omega\mid \diff{\mu}{\nu}(\omega)>y\})\,dy
\\ \dsp \phantom{m((0,\gamma])}
=\frac1{\mu(\Omega)}\int_{F_{\mu}(y)\ge (1-\gamma)\nu(\Omega)}
(F_{\mu}(+\infty)-F_{\mu}(y))\,dy,
\ \gamma\in(0,1].
\arre
\end{equation}
Since $\dsp F_{\mu}^{-1}$ is defined to be right continuous,
we have
\begin{equation}
\eqna{abscontigeneratedmeassprealsprf}
\arrb{l}\dsp
m((0,\gamma))
=\frac1{\mu(\Omega)}\int_{F_{\mu}(y)> (1-\gamma)\nu(\Omega)}
(F_{\mu}(+\infty)-F_{\mu}(y))\,dy
\\ \dsp \phantom{m((0,\gamma))}
=\frac1{\mu(\Omega)}\int_{y>F_{\mu}^{-1}((1-\gamma)\nu(\Omega))}
(F_{\mu}(+\infty)-F_{\mu}(y))\,dy
\\ \dsp \phantom{m((0,\gamma))}
=\frac1{\mu(\Omega)}\int_{F_{\mu}^{-1}((1-\gamma)\nu(\Omega))}^{\infty}
(F_{\mu}(+\infty)-F_{\mu}(y))\,dy
\arre
\end{equation}
for an open interval $(0,\gamma)$.

Fubini's Theorem implies,
\[\arrb{l}\dsp
m((0,1])=\frac1{\mu(\Omega)}\int_{\preals}
\nu(\{\omega\in\Omega\mid \diff{\mu}{\nu}(\omega)>y\})\,dy
=\frac1{\mu(\Omega)}\int_{\diff{\mu}{\nu}(\omega)>y}\,\nu(d\omega)\,dy
\\ \dsp \phantom{m((0,1])}
=\frac1{\mu(\Omega)}\int_{\Omega}
\biggl(\int_0^{\diff{\mu}{\nu}(\omega)}\,dy\biggr)\,\nu(d\omega)
=\frac1{\mu(\Omega)}\int_{\Omega}\diff{\mu}{\nu}(\omega)\,\nu(d\omega)
=\frac1{\mu(\Omega)}\int_{\Omega}\mu(d\omega)=1,
\arre\]
hence, $m$ is a probabiity measure on $(0,1]$.
\lemb
\lema{abscontigeneratedmeassprealsprf}
It follows that
\begin{equation}
\eqna{rskmSKzeta}
g(\gamma)
:=\int_{[1-\gamma,1]} \frac{m(d\alpha)}{\alpha}
=\frac{\nu(\Omega)}{\mu(\Omega)}\,F^{-1}_{\mu}(\nu(\Omega)\gamma),
\ \gamma\in[0,1).
\end{equation}
\DDD
\leme
\prfb
Not that for $\gamma\in[0,1)$ 
\[\arrb{l}\dsp
\{(\alpha,\beta)\in(0,1]\times[0,1)\mid
1-\gamma\le \beta<1,\ 1-\beta\le \alpha\le1\}
\\ \dsp {}
=\{(\alpha,\beta)\in(0,1]\times[0,1)\mid 0<\alpha\le 1,
\ 1-(\alpha\wedge\gamma) \le \beta<1\},
\arre\]
hence
$\dsp
\int_{[1-\gamma,1)}g(\beta)\,d\beta
=\int_{(0,1]}  \frac{\alpha\wedge\gamma}{\alpha}\,m(d\alpha)
= m((0,\gamma))+\gamma\,g(1-\gamma)$\,.

On the other hand, substituting $F=F_{\mu}$\,,
$\dsp \alpha=(1-\gamma)\nu(\Omega)$, and $\dsp a=F_{\mu}^{-1}(\alpha)$
in \eqnu{invfcnintglext}, and then using 
\eqnu{abscontigeneratedmeassprealsprf}.
and changing integration variables as $\dsp \beta=\beta'\nu(\Omega)$,
and also using $\dsp F_{\mu}(+\infty)=\nu(\Omega)$, we have
\[
\frac{\nu(\Omega)}{\mu(\Omega)}
\int_{1-\gamma}^{1} F_{\mu}^{-1}(\beta'\,\nu(\Omega))\,d\beta'
-m((0,\gamma))=\frac{\nu(\Omega)}{\mu(\Omega)}
\,\gamma\,F_{\mu}^{-1}((1-\gamma)\nu(\Omega)),
\ 0\le\gamma\le 1.
\]

Therefore, if we define $h:\ (0,1]\to\reals$ by
$\dsp h(\gamma)=\frac{\nu(\Omega)}{\mu(\Omega)}
\,F_{\mu}^{-1}((1-\gamma)\nu(\Omega))-g(1-\gamma)$,
we have
$\dsp
\frac1x\int_{(0,x]}h(\beta)\,d\beta=h(x),\ x\in(0,1],
$
which implies that $h$ is a constant.
Substituting $\gamma=1$ in \eqnu{rskmSKm} and 
\eqnu{abscontigeneratedmeassprealsprf},
we see that 
\[
g(0)=m(\{1\})
=\frac1{\mu(\Omega)}\int_{F_{\mu}(y)= 0}
(F_{\mu}(+\infty)-F_{\mu}(y))\,dy
=\frac{\nu(\Omega)}{\mu(\Omega)}\,F_{\mu}^{-1}(0)
\]
which implies $h(1)=0$. Therefore \eqnu{rskmSKzeta} holds.
\QED
\prfe
\prpb
Let $(\Omega,{\cal F})$ be a measurable space satisfying $\sgcn\ne\emptyset$,
and $\nu\in {\cal M}(\Omega)$.
Let $\mu\in {\cal M}(\Omega)$ be a finite measure satisfying $\mu\ll\nu$,
and $v_{\mu}:\ {\cal F}\to\reals$ be the submodular function 
\eqnu{fcnintgl2submodular},
and $m$ be the measure on $(0,1]$ defined in \eqnu{rskmSKm}.
Then, 
for a non-negative valued measurable function $f:\ \Omega\to\preals$,
the Choquet integration $v_{\mu}(f)$ of $f$ with respect to $v_{\mu}$ 
satisfies,
\begin{equation}
\eqna{rskmSKm2vf}
\arrb{l}\dsp
v_{\mu}(f)=\int_{(0,1]} v_{\alpha}(f)\, m(d\alpha),
\mbox{ where, }\\ \dsp 
v_{\alpha}(f)=
\frac{\mu(\Omega)}{\alpha\nu(\Omega)}
\int_0^{\infty} \nu(\{\omega\in\Omega\mid
f(\omega)>z\}) \wedge (\alpha\nu(\Omega))\,dz,\ 0<\alpha\le 1.
\arre
\end{equation}
\DDD
\prpe
\prfb
The definition \eqnu{fcnintgl2submodular} of $\dsp v_{\mu}$ and
the definition \eqnu{Choquetbd} of Choquet integration imply,
as in the proof of 
\eqnu{genDKChoquetsubmodular} in \thmu{uniformsubmodularfcn},
\begin{equation}
\eqna{vmuf}
v_{\mu}(f)
=\int_{0}^{\infty} 
\int_0^{\infty} (\nu(\Omega)-F_{\mu}(y))\wedge
\nu(\{\omega\in\Omega\mid f(\omega)>z\})\,dy\,dz.
\end{equation}

On the other hand, using
\eqnu{abscontigeneratedmeassprealsprf} with 
$\dsp \gamma=1-\frac{F_f(z)}{\nu(\Omega)}$
and
\eqnu{rskmSKzeta} with
$\dsp \gamma=\frac{F_f(z)}{\nu(\Omega)}$\, we have
\[\arrb{l}\dsp
\int_{(0,1]}\nu(\{\omega\in\Omega\mid f(\omega)>z\})
 \wedge (\alpha\nu(\Omega))
\frac{m(d\alpha)}{\alpha\nu(\Omega)}
\\ \dsp {}
=m((0,1-\frac{F_f(z)}{\nu(\Omega)}))
+(1-\frac{F_f(z)}{\nu(\Omega)})
\int_{[1-\frac{F_f(z)}{\nu(\Omega)},1]} \frac{m(d\alpha)}{\alpha}
\\ \dsp {}
=\frac1{\mu(\Omega)}\int_{F_{\mu}^{-1}(F_f(z))}^{\infty}
(\nu(\Omega)-F_{\mu}(y))\,dy
+\frac{1}{\mu(\Omega)}\,F^{-1}_{\mu}(F_f(z))\,
(\nu(\Omega)-F_f(z)).
\arre\]
Using the right continuity of $\dsp F_{\mu}^{-1}$ and the fact that
a point has zero $1$-dimensional Lebesgue measure
in a similar way as in the proof of \prpu{FubiniChoquetintegration},
we see that the right hand side is equal to
$\dsp
\frac1{\mu(\Omega)}\int_0^{\infty} (\nu(\Omega)-F_{\mu}(y))\wedge
\nu(\{\omega\in\Omega\mid f(\omega)>z\})\,dy
$.
Therefore 
\[\arrb{l}\dsp
\int_{(0,1]}\nu(\{\omega\in\Omega\mid f(\omega)>z\})
 \wedge (\alpha\nu(\Omega))
\frac{m(d\alpha)}{\alpha\nu(\Omega)}
\\ \dsp {}
=\frac1{\mu(\Omega)}\int_0^{\infty} (\nu(\Omega)-F_{\mu}(y))\wedge
\nu(\{\omega\in\Omega\mid f(\omega)>z\})\,dy.
\arre\]
Integrating over $z\ge0$ and using \eqnu{vmuf}, we have \eqnu{rskmSKm2vf}.
\QED
\prfe
%
%
We gave several formulas expressing values of submodular functions for sets
or Choquet integrations of functions in terms of supremum over class of 
measures, such as
\eqnu{submodextremeptrep} and
\eqnu{submodulartotalcorecoherentriskdef} of 
\thmu{totorder2modularvarrepbigcoreequiv} in \secu{submodformula},
\eqnu{genDelbaenKusuoka} of
\thmu{genDelbaenKusuoka} in \secu{Choquetint},
and
\eqnu{DelbaenKusuokaf2} of \thmu{comonotoneChoquetExprep} 
in \secu{lawinvcrm}.
Note that formulas in \secu{lawinvcrm} for 
$\nu$-uniform (`law invariant') submodular functions
and those in \secu{submodformula} and \secu{Choquetint} 
have different classes of measures to take supremum.
How far they are to be mathemtaically united is not clear to the author.
%

\section{Recursion relation associated with the fundamental formula}
\seca{submodextremeptreprecursion}

Let us return to the recursion relation \eqnu{submodextremeptreprecursion}
in \secu{intro} for set functions
on a measurable space $(\Omega,{\cal F})$ satisfying $\sgcn\ne\emptyset$.
\dfnu{vconti} of continuity of set function in this paper
implies that if $v_n$ is continuous, 
then the measures $\dsp \mu_{v_n,{\cal I}}$, ${\cal I}\in\sgcn$,
in the right hand side of \eqnu{submodextremeptreprecursion} 
are continuous ($\sigma$-additive) and well-defined as measures, 
hence \eqnu{submodextremeptreprecursion} makes sense.

We do not know in general if the resulting set function $v_{n+1}$ is
continous. On the other hand, if $v_n$\,, $n=0,1,2,\ldots$, are
all continuous for some $(\Omega,{\cal F})$ and $v_0$\,, then
the recursion converges.
\prpb
\prpa{recursionconv}
Let $(\Omega,{\cal F})$ be a meassurable space satisfying $\sgcn\ne\emptyset$,
and $\dsp v_0:\ {\cal F}\to\preals$ a non-decreasing and continuous 
set function satisfying $v_0(\emptyset)=0$.
If $\dsp v_n:\ {\cal F}\to\preals$, $n=0,1,2,\ldots$, is a sequence of
continuous set functions which satisfies 
\eqnu{submodextremeptreprecursion},
then for each $A\in{\cal F}$,
the sequence $\dsp v_n(A)$, $n=0,1,2,\ldots$, in non-decreasing
and bounded from above by $v_0(\Omega)$, hence converges. 
The limit $v:\ {\cal F}\to\preals$ defined by
$\dsp v(A)=\limf{n}v_n(A)$, $A\in{\cal F}$,
is non-decreasing, submodular, and satisfies \eqnu{upperconti} and
$v(\emptyset)=0$ and $v(\Omega)=v_0(\Omega)$.
\DDD
\prpe
\prfb
Since we assume continuity of $v_n$ for all $n\in\pintegers$,
$\dsp \mu_{v_n,{\cal I}}$ are measures.
The recursion relation therefore implies that $\dsp v_{n+1}$
is non-decreasing and 
$\dsp v_{n+1}(\emptyset)=0$ and $v_{n+1}(\Omega)=v_0(\Omega)$
for all $n\in\pintegers$, hence 
$\dsp v(\emptyset)=\limf{n}v_{n}(\emptyset)=0$ and
$\dsp v(\Omega)=\limf{n}v_{n}(\Omega)=v_0(\Omega)$.
Since $\dsp v_{n+1}$ is non-decrasing,
$\dsp v_{n+1}(A)\le v_{n+1}(\Omega)=v_0(\Omega)$ for all $n$ and $A$.
For $A\in{\cal F}$, Note that ${\cal I}_A\in\sgcn$, the insertion of $A$ to
${\cal I}\in\sgcn$ defined by \eqnu{insertion0}
satisfies $A\in{\cal I}_A$\,, which, with \eqnu{extremept}, implies
\[
v_{n+1}(A)=\sup_{{\cal I}\in \sgcn} \mu_{v_n,{\cal I}}(A)
\ge \mu_{v_n,{\cal I}_A}(A)=v_n(A),
\]
hence the sequence $\dsp v_n(A)$, $n=0,1,2,\ldots$, is non-decreasing.
We have seen that the sequence is bounded from above, hence it converges.
Thus 
\begin{equation}
\eqna{recursionconvprf2}
v(A)=\limf{n}v_n(A)=\sup_{n\in\nintegers} v_n(A),\ A\in{\cal F}.
\end{equation}
The non-decreasing property of measures and \eqnu{submodextremeptreprecursion}
imply that the limit $v$ is also non-decreasing.

To prove that $v$ is submodular, let $A,B\in{\cal F}$.
Then $\dsp v(A\cap B)=\limf{n}v_n(A\cap B)$ 
and $\dsp v(A\cup B)=\limf{n}v_n(A\cup B)$ imply that
for any $\eps>0$ there exists $n\in\nintegers$ such that
\begin{equation}
\eqna{recursionconvprf1}
v(A\cap B)+v(A\cup B)\le v_{n}(A\cap B)+v_n(A\cup B)+\eps.
\end{equation}
Let $\dsp ({\cal I}_A)_B\in \sgcn$ be the sequential insertion of 
$A$ and $B$ to ${\cal I}$ in \eqnu{insertionAB}.
Then as remarked below \eqnu{insertionAB}, 
$\dsp A\cap B\in({\cal I}_A)_B$ and
$\dsp A\cup B\in({\cal I}_A)_B$\,,
hence by definition of extremal measure in \eqnu{extremept},
additivity of measure, \eqnu{submodextremeptrep},
and \eqnu{recursionconvprf2}
it holds that
\[\arrb{l}\dsp
v_n(A\cap B)+v_n(A\cup B)
=\mu_{v_n,({\cal I}_A)_B}(A\cap B)+\mu_{v_n,({\cal I}_A)_B}(A\cup B)
\\ \dsp \phantom{v_n(A\cap B)+v_n(A\cup B)}
=\mu_{v_n,({\cal I}_A)_B}(A)+\mu_{v_n,({\cal I}_A)_B}(B)
\\ \dsp \phantom{v_n(A\cap B)+v_n(A\cup B)}
\le v_{n+1}(A)+v_{n+1}(B)
\\ \dsp \phantom{v_n(A\cap B)+v_n(A\cup B)}
\le v(A)+v(B).
\arre\]
Since $\eps>0$ is arbitrary,
this and \eqnu{recursionconvprf1} proves the submodularity
$\dsp v(A\cap B)+v(A\cup B)\le v(A)+v(B)$.

To prove \eqnu{upperconti}
asuume that ${\cal A}=\{A_n\mid n=1,2,3,\ldots\}\subset{\cal F}$ satisfies
$A_1\subset A_2\subset\cdots$.
Since by assumption $v_k$ is continuous, and since we saw that $v_k$ is
non-decreasing, \prpu{conti24toconti23} implies 
\eqnu{upperconti} for $v_k$:
$\dsp
\limf{n} v_k(A_n)=v_k(\bigcup_{n\in\nintegers} A_n),\ k\in\pintegers.
$
This and $\dsp v(\bigcup_{n\in\nintegers} A_n)=
\limf{k} v_k(\bigcup_{n\in\nintegers} A_n)$ imply that
for any $\eps>0$ there exists $k\in\pintegers$ such that
\[
v(\bigcup_{n\in\nintegers} A_n)
\le v_k(\bigcup_{n\in\nintegers} A_n)+\eps
= \limf{n} v_k(A_n)+\eps.
\]
Since $v$ is non-decreasing and
$\dsp v(A_n)=\sup_{k\in\nintegers} v_k(A_n)$ we
further have
\[
\limf{n} v(A_n)\le v(\bigcup_{n\in\nintegers} A_n)
\le \limf{n} v(A_n)+\eps.
\]
Since $\eps>0$ is arbitraly this proves \eqnu{upperconti}.
\QED
\prfe
Note that we assume continuity of each $v_n$ in this result.
The continuity of the limit $v$, or equivalently, \eqnu{lowerconti} for $v$,
is an open problem.

If the total set is a finite set, then
there are only finite number of distinct subsets, hence
if a set function is finitely additive then it is $\sigma$-additive
(continuous). In other words, any finitely additive measures is a
($\sigma$-additive) measure if the total set is a finite set. Therefore, 
\dfnu{vconti} of continuity of set function
is always satisfied, and in particular,
the recursion relation \eqnu{submodextremeptreprecursion} makes sense
for any non-decreasing initial set function $v_0$ and all $n\in\nintegers$.

Hereafter, we assume that the cardinality $m$ of the total set 
$\dsp\Omega=\Omega_m$ is finite: $\dsp m\in\nintegers$,
and assume for simplicity $\dsp{\cal F}=2^{\Omega_m}$, and
consider the recursion relation
\eqnu{submodextremeptreprecursion}:
$\dsp v_{n+1}=\sup_{{\cal I}\in \sgcn} \mu_{v_n,{\cal I}}$\,,
$n=0,1,2,\ldots$. with $\dsp v_0:\ 2^{\Omega_m}\to\preals$ a
non-decreasing (automatically continuous) set function satisfying
$\dsp v_0(\emptyset)=0$.
\prpu{recursionconv} implies that $v_n(A)$, $n\in\pintegers$,
is non-decreasing in $n$ and converges as $n\to\infty$, 
for each $A\in{\cal F}$.

Note that if we put $\dsp \Omega_m=\{\omega_i\mid i=1,\ldots,m\}$, then 
${\cal I}\in\sgcn$ is of a form
\begin{equation}
\eqna{submodextremeptreprecursion1endprf1}
{\cal I}=\{\emptyset,\{\omega_{i_1}\},\{\omega_{i_1},\omega_{i_2}\},
\ldots,\{\omega_{i_1},\ldots,\omega_{i_{m-1}}\},\Omega_m\}.
\end{equation}
In particular, the cardinality $\dsp \sharp\sgcn$ of $\sgcn$ satisfies
$\dsp\sharp\sgcn\le m!$. 
There are at most $m!$ numbers to compare in the supremum
in the right hand side of \eqnu{submodextremeptreprecursion},
hence supremum is always attained; for each $A\subset \Omega_m$
there exists ${\cal I}\in\sgcn$ such that
\begin{equation}
\eqna{submodextremeptreprecursion1endprf2}
v_{n+1}(A)= \mu_{v_n,{\cal I}}(A).
\end{equation}

If the limit is attained after a finite iteration of
\eqnu{submodextremeptreprecursion},
namely, if $\dsp v_{n+1}=v_n$ for an integer $n$,
then \eqnu{submodextremeptreprecursion} implies
$\dsp v_{n}=\sup_{{\cal I}\in \sgcn} \mu_{v_n,{\cal I}}$\,,
which implies \eqnu{submodextremeptrep} for $\dsp v=v_n$
and \thmu{totorder2modularvarrepbigcoreequiv} then implies that
$v_n$ is submodular.
On the other hand, if we start from $v_0$ which is not submodular,
then \thmu{totorder2modularvarrepbigcoreequiv} implies $v_1\ne v_0$\,.
Thus the recursion relation \eqnu{submodextremeptreprecursion} is of
interest in relating non-submodular functions to submodular functions.

As a first elementary example, 
we will prove that if $m=3$, then for any non-decreasing $v_0$
we have $v_2=v_1$\,, so that the limit is reached in at most $1$ step 
recursion and $v_1$ is submodular.
\prpb
\prpa{submodextremeptreprecursion1end}
Let $m\in\nintegers$ and $\Omega_m$ be a set of cardinality 
$\dsp \sharp\Omega_m=m$.
\itmb
\item
For any set function $\dsp v:\ 2^{\Omega_m}\to\preals$
if one of the following (i)--(vii)
 holds for a pair of sets $A,B\subset \Omega_m$,
then $\dsp v(A\cap B)+v(A\cup B)=v(A)+v(B)$ holds;
(i) $A=\emptyset$, or (ii) $B=\emptyset$, or (iii) $A=\Omega_m$, 
or (iv) $B=\Omega_m$, or (v) $A\subset B$, or (vi) $A\supset B$, or 
(vii) $A=B$.

\item
Let $\dsp v_0:\ 2^{\Omega_m}\to\preals$ be a non-decreasing set function
satisfying $\dsp v_0(\emptyset)=0$, and for $n\in\nintegers$,
let $v_n:\ 2^{\Omega_m}\to\preals$ be the non-decreasing set function
determined by the recursion relation \eqnu{submodextremeptreprecursion}.
Then for any positive integer $n\in\nintegers$,
if one of the following (i)--(iv) holds for a pair of sets $A,B\subset \Omega_m$,
then $\dsp v_n(A\cap B)+v_n(A\cup B)\le v_n(A)+v_n(B)$ holds;
(i) $A\cap B=\emptyset$,
(ii) $A\cup B=\Omega_m$,
(iii) $\sharp A=1$,
(iv) $\sharp A=m-1$,
where 
$\sharp A$ denotes the cardinality of (number of elements in) the set $A$.
\DDD
\itme
\prpe
\prfb
\itmb
\item
By interchanging $A$ and $B$, it suffices to prove (i), (iv), (v).
If $A=\emptyset$ then 
$\dsp v(A\cap B)+v(A\cup B)=v(\emptyset)+v(B)=v(A)+v(B)$,
If $B=\Omega_m$ then 
$\dsp v(A\cap B)+v(A\cup B)=v(A)+v(\Omega_m)=v(A)+v(B)$,
If $A\subset B$ then $\dsp v(A\cap B)+v(A\cup B)=v(A)+v(B)$.

\item
(i) If $A\cap B=\emptyset$, then
as noted at \eqnu{submodextremeptreprecursion1endprf2},
there exists ${\cal I}\in \sgcn$ such that
$\dsp v_{n}(A\cup B)= \mu_{v_{n-1},{\cal I}}(A\cup B)$.
Using the additivity of the measure $\dsp \mu_{v_{n-1},{\cal I}}$ with
the assumption $A\cap B=\emptyset$, we have
\[\arrb{l}\dsp
 v_n(A\cap B)+v_n(A\cup B)= \mu_{v_{n-1},{\cal I}}(A\cup B)
= \mu_{v_{n-1},{\cal I}}(A)+ \mu_{v_{n-1},{\cal I}}(B)
\\ \dsp {}
\le v_n(A)+v_n(B).
\arre\]

(ii) If $A\cup B=\Omega_m$, then as in the case (i)
there exists ${\cal I}\in \sgcn$ such that
$\dsp v_{n}(A\cap B)= \mu_{v_{n-1},{\cal I}}(A\cap B)$.
Noting that $\dsp v_n(\Omega_m)=v_0(\Omega_m)=\mu_{v_{n-1},{\cal I}}(\Omega_m)$
for all ${\cal I}\in \sgcn$, and
using the aditivity of the measure $\dsp \mu_{v_{n-1},{\cal I}}$, we have
\[\arrb{l}
 v_n(A\cap B)+v_n(A\cup B)= \mu_{v_{n-1},{\cal I}}(A\cap B)+
\mu_{v_{n-1},{\cal I}}(\Omega_m)
\\ \dsp {}
= \mu_{v_{n-1},{\cal I}}(A)+ \mu_{v_{n-1},{\cal I}}(B)
\le v_n(A)+v_n(B).
\arre
\]

(iii) If $\sharp A=1$ then,
either $A\cap B=\emptyset$ or $A\subset B$ holds for any $B\subset \Omega_m$.
In both cases we already proved that
$\dsp v_n(A\cap B)+v_n(A\cup B)\le v_n(A)+v_n(B)$ holds.

(iv) If $\sharp A=m-1$ then,
either $A\cup B=\Omega_m$ or $A\supset B$ holds for any $B\subset \Omega_m$.
In both cases we already proved that
$\dsp v_n(A\cap B)+v_n(A\cup B)\le v_n(A)+v_n(B)$ holds.
\QED
\itme
\prfe
\corb
\cora{submodextremeptreprecursion1end}
If $m=\sharp\Omega_m=3$ and $v_0:\ 2^{\Omega_3}\to\preals$, 
is non-decreasing and $\dsp v_0(\emptyset)=0$,
then the sequence $v_n:\ 2^{\Omega_3}\to\preals$, $n\in\nintegers$, 
determined by the recursion relation \eqnu{submodextremeptreprecursion}
satisfies $v_1=v_2=\cdots$ and $v_1$ is a submodular function.
\DDD
\core
\prfb
If $A\subset\Omega_3$, then $\dsp 0\le \sharp A\le \sharp \Omega_3=3$.
This implies $A=\emptyset$ or $A=\Omega_3$ or $\sharp A=1$ or $\sharp A=2=m-1$.
\prpu{submodextremeptreprecursion1end} implies that for all the cases,
$\dsp v_n(A\cap B)+v_n(A\cup B)\le v_n(A)+v_n(B)$ holds
for any $B\subset \Omega_3$ and for all $n=1,2,\ldots$.
Therefore $\dsp v_1$ is submodular and
the recursion relation reaches the limit at $n=1$.
\QED
\prfe
It turns out that not only for $m=3$, but for all $m\in\nintegers$
we can easily find examples that
the recursion relation \eqnu{submodextremeptreprecursion} reaches
the limit after just $1$ iteration.
\prpb
\prpa{submodextremeptreprecursion1endfiniteuniform}
Let $\dsp m=\sharp\Omega_m\in\nintegers$ and
$\dsp v_0:\ 2^{\Omega_m}\to\preals$ be non-decreasing and
$\dsp v_0(\emptyset)=0$, and assume that there exists 
$\dsp g:\ \{1,2,\ldots,m\}\to\preals$ such that
$\dsp v_0(A)=g(\sharp A)$, $A\subset \Omega_m$.
Then the recursion \eqnu{submodextremeptreprecursion} reaches the limit
at $n=1$, i.e., $\dsp v_2=v_1$, and hence $v_1$ is submodular.
\DDD
\prpe
\prfb
Let $S_m$ denote the collection of the rearrangements of the 
sequence of $m$ elements of the total set $\Omega_m$;
\begin{equation}
\eqna{rearragement}
S_m=\{(\sigma_1,\ldots,\sigma_m)\mid
\{\sigma_1,\ldots,\sigma_m\}=\Omega_m\},
\end{equation}
and for $\dsp \sigma=(\sigma_1,\ldots,\sigma_m)\in S_m$ put
\begin{equation}
\eqna{rearrangementI}
{\cal I}_{\sigma}=\{\emptyset,\{\sigma_{1}\},\{\sigma_{1},\sigma_{2}\},
\ldots,\{\sigma_{1},\ldots,\sigma_{m-1}\},\Omega_m\}.
\end{equation}
As discussed at around \eqnu{submodextremeptreprecursion1endprf1} we have
\begin{equation}
\eqna{sgcnfiniteset}
\sgcn=\{{\cal I}_{\sigma}\mid \sigma\in S_m\}.
\end{equation}

The definition \eqnu{extremept} of extremal measure implies,
with the assumption $\dsp v_0(A)=g(\sharp A)$,
\begin{equation}
\eqna{terminalmeasfiniteset}
\mu_{v_0,{\cal I}_{\sigma}}(\{\sigma_i\})
=v_0(\{\sigma_1,\ldots,\sigma_i\})-v_0(\{\sigma_1,\ldots,\sigma_{i-1}\})
=g(i)-g(i-1),\  \ i=1,\ldots,m.
\end{equation}
Note that $\dsp g(0)=v_0(\emptyset)=0$.
Since $v_0$ is non-decreasing, so is $g$.
In particular, $g(i)-g(i-1)\ge0$, $i=1,\ldots,m$.

Choose the rearrangement $\alpha=(\alpha_1,\ldots,\alpha_m)$ of
first $m$ positive integers $\{1,\ldots,m\}$ such that
\begin{equation}
\eqna{submodextremeptreprecursion1endfiniteuniformprf1}
g(\alpha_1)-g(\alpha_1-1)\ge g(\alpha_2)-g(\alpha_2-1)\ge\cdots
\ge g(\alpha_m)-g(\alpha_m-1)\ge0.
\end{equation}
In particular, \eqnu{terminalmeasfiniteset} and
\eqnu{submodextremeptreprecursion1endfiniteuniformprf1} imply
$\dsp g(\alpha_j)-g(\alpha_j-1)
=\mu_{v_0,{\cal I}_{\sigma}}(\{\sigma_{\alpha_j}\})$, $\sigma\in S_m$\,.

Put $\dsp \Omega_m=\{\omega_1,\ldots,\omega_m\}$.
Since $\dsp\{\sigma_{\alpha_j}\mid j=1,\ldots,m\}=\Omega_m$,
there exists $\sigma\in S_m$ dependent rearrangement of
first $m$ positive integers
$\dsp h_{\sigma}:\ \{1,\ldots,m\}\to\{1,\ldots,m\}$ such that
$\dsp \sigma_{\alpha_{h_{\sigma}(i)}}=\omega_i$, $i=1,\ldots,m$.
Then
\begin{equation}
\eqna{submodextremeptreprecursion1endfiniteuniformprf3}
\mu_{v_0,{\cal I}_{\sigma}}(\{\omega_i\})=
g(\alpha_{h_{\sigma}(i)})-g(\alpha_{h_{\sigma}(i)}-1),\ i=1,\ldots,m.
\end{equation}
Hence, for $\dsp A\in2^{\Omega_m}$,
\[
v_1(A)=\sup_{\sigma\in S_m} \mu_{v_0,{\cal I}_{\sigma}}(A)
=\sup_{\sigma\in S_m} \sum_{\omega_i\in A}
 \mu_{v_0,{\cal I}_{\sigma}}(\{\omega_i\})
=\sup_{\sigma\in S_m} \sum_{\omega_i\in A} 
(g(\alpha_{h_{\sigma}(i)})-g(\alpha_{h_{\sigma}(i)}-1)).
\]
Supremum is attained by choosing largest $\sharp A$ terms in
\eqnu{submodextremeptreprecursion1endfiniteuniformprf1}, hence
\begin{equation}
\eqna{submodextremeptreprecursion1endfiniteuniformprf2}
v_1(A)= \sum_{j=1}^{\sharp A}(g(\alpha_{j})-g(\alpha_{j}-1)).
\end{equation}

We can repeat the calculations by substituting $v_0$ for $v_1$, to find
\[
\mu_{v_1,{\cal I}_{\sigma}}(\{\omega_i\})=
g(\alpha_{h_{\sigma}(i)})-g(\alpha_{h_{\sigma}(i)}-1),\ i=1,\ldots,m,
\]
whose right hand side is equal to 
\eqnu{submodextremeptreprecursion1endfiniteuniformprf3}, hence
\[
v_2(A)=\sup_{\sigma\in S_m} \mu_{v_1,{\cal I}_{\sigma}}(A)
= \sum_{j=1}^{\sharp A}(g(\alpha_{j})-g(\alpha_{j}-1))
=v_1(A),\ A\in{\cal F},
\]
which implies $\dsp v_1=\sup_{{\cal I}\in\sgcn} \mu_{v_1,{\cal I}}$.
\QED
\prfe

With all these results for examples of $v_2= v_1$\,,
it may be interesting to see that there is an example of $v_2\ne v_1$
for $m=4$. The table below gives such an example.
For example, we see from the table that 
$\dsp v_1(\{1,2\})=17\ne18=v_2(\{1,2\})$, hence $v_1$ is not submodular
and it is not the limit of the recursion.
We can also see 
\[
v_1(\{1\})+v_1(\{1,2,3\})=11+31=42>41=17+24=v_1(\{1,2\})+v_1(\{1,3\})
\]
which directly proves that $v_1$ is not submodular.
\par\noindent
\tabcolsep 1pt \tabskip 0pt
\tabb{|c||c|c|c|c|c|c|c|c|c|c|c|c|c|c|c|c|}\hline
A&\,$\emptyset$\,&$\{1\}$&$\{2\}$&$\{3\}$&$\{4\}$&$\{1,2\}$&$\{1,3\}$&$\{1,4\}$&$\{2,3\}$&$\{2,4\}$&$\{3,4\}$&$\{1,2,3\}$&$\{1,2,4\}$&$\{1,3,4\}$&$\{2,3,4\}$&$\Omega_4$\\
\hline\hline
$v_0(A)$&$ 0$&$ 7$&$13$&$20$&$19$&$17$&$24$&$30$&$28$&$34$&$41$&$31$&$36$&$42$&$43$&$44$\\\hline
$v_1(A)$&$ 0$&{\bf $11$}&{\bf $15$}&{\bf $22$}&{\bf $23$}&$17$&$24$&$30$&$28$&$34$&$41$&$31$&$36$&$42$&$43$&$44$\\\hline
$v_2(A)$&$ 0$&$11$&$15$&$22$&$23$&{\bf $18$}&{\bf $25$}&$30$&{\bf $29$}&$34$&$41$&$31$&$36$&$42$&$43$&$44$\\\hline
\tabe
\par\noindent
Incidentally, in this example we have $v_3=v_2$\,, so that 
the recursion \eqnu{submodextremeptreprecursion} reaches the limit
at $n=2$, and $v_2$ is submodular.
It is an open problem whether there exists an infinite sequence $v_n$\,,
$n=0,1,2,\ldots$, of non-decreasing and non-submodular sset functions
which satisfies the recursion \eqnu{submodextremeptreprecursion}.

\end{document}